\documentclass{amsart}

\usepackage{fullpage}
\usepackage[utf8]{inputenc}
\usepackage[T1]{fontenc}
\usepackage{graphicx}
\usepackage{amsmath,amsfonts,amssymb,mathtools}
\usepackage{stmaryrd}
\usepackage{graphicx}
\usepackage{mathrsfs,dsfont}

\usepackage{hyperref}

\usepackage{enumerate}

\usepackage{color}

\newtheorem{theo}{Theorem}[section]
\newtheorem{algo}[theo]{Algorithm}
\newtheorem{cor}[theo]{Corollary}
\newtheorem{rem}[theo]{Remark}

\newtheorem{propo}[theo]{Proposition}
\newtheorem{lemme}[theo]{Lemma}
\newtheorem{defi}[theo]{Definition}

\newtheorem{hyp}[theo]{Assumption}
\newcommand{\E}{\mathbb{E}}
\newcommand{\R}{\mathbb{R}}
\newcommand{\PP}{\mathbb{P}}

\newcommand{\N}{\mathbb{N}}

\begin{document}

\title{Analysis of Adaptive Multilevel Splitting algorithms in an idealized case}
\author{Charles-Edouard Br\'ehier}\address{Universit\'e Paris-Est, CERMICS (ENPC), INRIA, 6-8-10 avenue Blaise Pascal, 77455 Marne-la-Vall\'ee, France}
\address{INRIA Paris-Rocquencourt, Domaine de Voluceau - Rocquencourt, B.P. 105 - 78153 Le Chesnay, France}
\email{brehierc@cermics.enpc.fr}
\email{lelievre@cermics.enpc.fr}
\email{mathias.rousset@inria.fr}
\author{Tony Leli\`evre}
\author{Mathias Rousset}
\date{}
\begin{abstract}
The Adaptive Multilevel Splitting algorithm \cite{cerou-guyader-07a} is a very powerful and versatile method to estimate rare events probabilities. It is an iterative procedure on an interacting particle system, where at each step, the $k$ less well-adapted particles among $n$ are killed while $k$ new better adapted particles are resampled according to a conditional law. We analyze the algorithm in the idealized setting of an exact resampling and prove that the estimator of the rare event probability is unbiased whatever $k$. We also obtain a precise asymptotic expansion for the variance of the estimator and the cost of the algorithm in the large $n$ limit, for a fixed $k$.
\end{abstract}
\maketitle
\tableofcontents
\section{Introduction}

Let $X$ be a real random variable, such that $X > 0$ almost surely. We want to approximate the following probability:
\begin{equation}
p=\PP(X \ge a),
\end{equation}
where $a>0$ is a given threshold, such that $p > 0$. When $a$ goes to infinity, the above probability goes to~$0$, meaning that $\left\{X \ge a\right\}$ becomes a rare event. Such problems appear in many contexts, such as molecular dynamics simulations~\cite{cerou-guyader-lelievre-pommier-11} or reliability problems with many industrial applications,  for example.

Estimating a rare event using a direct Monte Carlo estimation is inefficient, as can be seen by the analysis of the relative error. Indeed, let $(X_i)_{i\in\N}$ be a sequence of independent and identically distributed random variables with the same law as $X$. Then for any positive integer $M$,
\begin{equation}\label{estim_pureMC}
\hat{p}_M=\frac{1}{M}\sum_{n=1}^{M}\mathds{1}_{X_n \ge a}
\end{equation}
is an unbiased estimator of $p$: $\E[\hat{p}_M]=p$. It is also well-known that its variance is given by $\text{Var}(\hat{p}_M)=\frac{p(1-p)}{M}$ and therefore the relative error writes:
\begin{equation}\label{err_rel_pureMC}
\frac{\sqrt{\text{Var}(\hat{p}_M)}}{p}= \sqrt{\frac{1-p}{Mp}}.
\end{equation}
Assume that the simulation of one random variable $X_n$ requires a computational cost $c_0$. For a relative error of size $\epsilon$, the cost of a direct Monte Carlo method is thus of the order
\begin{equation}\label{cout_pureMC}
c_0\frac{1-p}{\epsilon^2 p}.
\end{equation}
which is prohibitive for small probabilities (say $10^{-9}$).

Many algorithms devoted to the estimation of the probability of rare events have been proposed. Here we focus on the so-called Adaptive Multilevel Splitting (AMS) algorithm~\cite{cerou-guyader-07a,cerou-guyader-07b}.

Let us explain one way to understand this algorithm (see~\cite{brehier-gazeau-goudenege-lelievre-rousset-2014} for a more general presentation). The splitting strategy relies on the following remark. Let us introduce $J$ intermediate levels: $a_0=0<a_1<\ldots <a_J=a$. The small probability $p$ satisfies:
$$p=\prod_{j=1}^J p_j$$
where $p_j = \PP(X>a_j | X>a_{j-1})$, with $1\leq j\leq J-1$ and $p_J=\PP(X\ge a | X>a_{J-1})$. In order to use this identity to build an estimator of $p$, one needs to (i) define appropriately the intermediate levels and (ii) find a way to sample according to the conditional distributions $\mathcal{L}(X | X>a_{j-1})$ to approximate each $p_j$ using independent Monte Carlo procedures.


In this article, we will be interesting in the idealized case where we assume we have a way to draw independent samples according to the conditional distributions $\mathcal{L}(X | X>a_{j})$, $j \in \{1, \ldots, J-1\}$. We will discuss at length this assumption below. It is then easy to check that for a given $J$, the variance is minimized when $p_1=\ldots=p_J=p^{1/J}$, and that the associated variance is a decreasing function of $J$. It is thus natural to try to devise a method to find the levels $a_j$ such that $p_1=\ldots=p_J$. In AMS algorithms, the levels are defined in an adaptive and random way in order to satisfy (up to statistical fluctuations) the equality of the factors $p_j$. This is based on an interacting particle system approximating the quantiles $\PP(X\ge a )$ using an empirical distribution. The version of the algorithm we study depends on two parameters: $n$ and $k$. The first one denotes the total number of particles. The second one denotes the number of resampled particles at each iteration: they are those with the $k$ lowest current levels (which mean that a sorting procedure is required). Thus,  the levels are defined in such a way that $p_j=\left(1-\frac k n \right)$ and the estimator of the probability $p$ writes:
$$\hat{p}_{n,k}= C^{n,k} \left( 1- \frac{k}{n}\right)^{J^{n,k}}$$
where ${J^{n,k}}$ is the number of iterations required to reach the target level $a$, and $C^{n,k} \in \left[1-\frac{k-1}{n},1\right]$ is a correction factor precisely defined below (see~\eqref{eq:corrector}). Notice that $C^{n,k}=1$ if $k=1$.

In all the following, we will make the following assumption:
\begin{hyp}\label{hyp:static}
$X$ is a real-valued positive random variable which admits a continuous cumulative distribution function $t \mapsto \PP(X \le t)$.
\end{hyp}
 This ensures for example that (almost surely), there is no more than one particle at the same level, and thus that the resampling step in the algorithm is well defined. We will show in Section~\ref{sec:hyp} below that this assumption can be relaxed: the continuity of the cumulative distribution function is actually only required on $[0,a)$, in which case the AMS algorithm still yields an estimate of $\PP(X \ge a)$ (with a large inequality). From Section~\ref{sec:expo_wp}, we will always work under  Assumption~\ref{hyp:static}, and thus we will always use for simplicity strict rather than large inequalities on $X$ (notice that under Assumption~\ref{hyp:static}, $p=\PP(X\ge a)=\PP(X > a)$).

Our aim in this article is twofold. First, we show that for any values of $n$ and $k$, the estimator $\hat{p}_{n,k}$ is an unbiased estimator of $p$ (see Theorem \ref{th:unbiased}):
\begin{equation}\label{eq:bias_intro}
\boxed{
\E(\hat{p}_{n,k})=p.
}
\end{equation}
Second, we are able to obtain an explicit asymptotic expression for the variance of the estimator $\hat{p}_{n,k}$ in the limit of large $n$, for fixed $k$ and $p$, and thus for the relative error (see Proposition~\ref{propo:var}, Equation~\eqref{eq:var}):
\begin{equation}\label{eq:var_intro}
\boxed{
\frac{\sqrt{\mathrm{Var}(\hat{p}^{n,k})}}{p}=\frac{-\log p}{\sqrt{n}}\left(1+\frac{\left(1 -\log (p)\right)(k-1)}{2n}+{\rm o}\left(\frac{1}{n}\right)\right)^{1/2}.
}
\end{equation}
Thus, if we consider the cost associated to the Monte Carlo estimator based on $M$ independent realizations of the algorithm, and if $M$ is chosen in such a way that the relative error is of order $\epsilon$, one ends up with the following asymptotic cost for the AMS algorithm (see Theorem~\ref{theo:cost}, Equations~\eqref{eq:costAMS} and~\eqref{eq:cost_asympt}): for fixed $k$ and $p$, in the limit of large $n$
\begin{equation}\label{eq:cost_intro}
\boxed{
\frac{c_0 + c_1\log(n)}{\epsilon^2}\left[ \left( \left(\log (p) \right)^2-\log(p)\right) +\frac{1}{n}\left(-\log (p)\left(k-1 \right)+\frac{1}{2}\left(\log (p)\right)^2-\frac{1}{2}\left(\log (p)\right)^3\right) +{\rm o}\left(\frac{1}{n}\right) \right]
}
\end{equation}
where $c_0$ denotes the cost for drawing one sample according to the conditional distributions $\mathcal{L}(X | X>x)$ (assumed to be independent of $x$, to simplify), and $c_1\log(n)$ is the cost associated with the sorting procedures involved in the algorithm. Here again the Landau symbol ${\rm o}$ in the above depend non-uniformly on $k$ and $p$. The two results~\eqref{eq:var_intro} and~\eqref{eq:cost_intro} should be compared to the corresponding formulae for direct Monte Carlo simulation~\eqref{err_rel_pureMC} and~\eqref{cout_pureMC} above. From these results, we conclude that, in this asymptotic regime:
\begin{enumerate}[(i)]
\item the choice $k=1$ is the optimal one in terms of variance and cost
\item AMS yields better result than direct Monte Carlo if
$$\left( 1 + \frac{c_1}{c_0}\log(n)\right) \left( \left(\log p \right)^2-\log p\right) < \frac{1-p}{p}$$
which will be the case for sufficiently small probability $p$.
\end{enumerate}

The assumption that we are able to sample according to the conditional distributions $\mathcal{L}(X | X>x)$ (idealized setting) is a severe limitation, from a practical viewpoint. Let us make three comments about this idealized setting. First, to the best of our knowledge, all the theoretical results which have been obtained so far in the literature~\cite{cerou-del-moral-furon-guyader-12,cerou-guyader-07a,cerou-guyader-07b,guyader-hengartner-matzner-lober-11} rely on such an assumption. It is believed that the qualitative conclusions obtained under this assumption are still meaningful for the actual algorithm used in practice, where the conditional distributions $\mathcal{L}(X | X>a_{j-1})$ are only approximately sampled (using for example Metropolis Hastings procedures). In some sense, in the idealized setting, one studies the optimal performance one could reach with this algorithm. Second, we will describe in Section~\ref{sect:samp} situations where this assumption makes sense in practice: this is in particular the case in the 1d dynamic setting described in Section~\ref{sec:1ddyn}. Third, in a paper in preparation~\cite{brehier-gazeau-goudenege-lelievre-rousset-2014}, we will actually show that it is possible to obtain an unbiased estimator of $p$ using AMS in a very general setting. In other words, the idealized setting is crucial to obtain estimates on the variances and the costs of the algorithm, but not to prove the unbiasedness property for the estimator of the rare event probability.

Let us now review the results known from the literature on AMS. As mentioned above, the AMS algorithm has been introduced in~\cite{cerou-guyader-07a,cerou-guyader-07b}, where it is proven that their estimator $\hat{p}_{n,k}$ indeed converges almost surely to $p$ (in the limit $n \to \infty$). A central limit theorem is also provided, in the asymptotic regime  where $\frac{k}{n}=p_0$ is fixed. In~\cite{guyader-hengartner-matzner-lober-11}, the authors consider the case $k=1$, prove the unbiasedness of the estimator and analyze the variance. In~\cite{cerou-del-moral-furon-guyader-12}, the authors analyze the large $n$ limit, with fixed $\frac{k}{n}=p_0$. In summary, our results differ from what has been proven before in two ways: we prove that $\hat{p}_{n,k}$ is an unbiased estimator for any $k \ge 1$ (and not only $k=1$) and we analyze the variance and the cost at a fixed $k$ (in the limit $n \to \infty$) and show that in this regime, $k=1$ is optimal. 

In addition to the new results presented in this paper, we would like to stress that the techniques of proof we use seem to be original. The main idea is to consider (under Assumption \ref{hyp:static}) the family $(P(x))_{x\in[0,a]}$ of conditional probabilities:
\begin{equation}
P(x)=\PP(X>a | X>x),
\end{equation}
and to define associated estimators $\hat{p}_{n,k}(x)$ thanks to the AMS algorithm. We are then able to derive an {\em{explicit functional equation on $\E[\hat{p}_{n,k}(x)]$ as a function of $x$}}; which follows from the existence of an explicit expression for the distribution of order statistics of independent random variables. We can then check the unbiased property (see Theorem \ref{th:unbiased})
$$\E[\hat{p}_{n,k}(x)]=P(x).$$
Finally noting that $P(0)=p$, we get~\eqref{eq:bias_intro}. To analyze the computational cost (Theorem~\ref{theo:cost}), we follow a similar strategy, with several more technical steps. First, we derive functional equations for the variance and the mean number of iterations in the algorithm. In general we are not able to give explicit expressions to the solutions of these equations, and we require the following auxiliary arguments:
\begin{enumerate}[(i)]
\item We show how one can relate the general case to the so-called exponential case, when $X$ has an exponential distribution with mean $1$;
\item We then prove that the solutions of the functional equations are solutions of linear Ordinary Differential Equations of order $k$;
\item We finally get asymptotic results on the solutions to these ODEs in the limit $n\rightarrow +\infty$, with fixed values of $p$ and $k$.
\end{enumerate}




The paper is organized as follows. In Section~\ref{sec:AMS}, we introduce the AMS algorithm and discuss the sampling of the conditional distributions ${\mathcal L}(X | X > z)$. In Section~\ref{sec:expo_wp}, we show how to relate the case of a general distribution for $X$ to the case when $X$ is exponentially distributed. We then prove that the algorithm is well defined, in the sense that it terminates in a finite number of iterations (almost surely). In Section \ref{sect:unbiased}, we show one of our main result, namely Theorem~\ref{th:unbiased} which states that the estimators $\hat{p}_{n,k}(x)$ are unbiased. Finally, in Section \ref{sect:cost}, we study the cost of the algorithm, with asymptotic expansions in the regime where $p$ and $k$ are fixed and $n$ goes to infinity. The proofs of the results of Section \ref{sect:cost} are postponed to Section \ref{sect:details}.

\section{The Adaptive Multilevel Splitting algorithm}
\label{sec:AMS}

After presenting the AMS algorithm in Section~\ref{sect:algo}, we discuss the fundamental assumption that we know how to sample according to the conditional distributions $\mathcal{L}(X | X>z)$ in Section \ref{sect:samp}. We will in particular show that this assumption is actually practical at least in one setting: the one-dimensional dynamic setting presented in Section~\ref{sec:1ddyn}. Finally, we discuss Assumption~\ref{hyp:static} in Section \ref{sec:hyp}  and show that it can be replaced by a less stringent hypothesis. This is particularly useful in the framework of the high-dimensional dynamic setting described in Section~\ref{sec:nddyn}.

We would like to stress that Sections $2.2$ and $2.3$ discuss practical 
aspects of the implementation of AMS and can be skipped if the reader is only interested in the two 
main results about the unbiasedness and cost of AMS, presented in Sections \ref{sect:unbiased} and \ref{sect:cost}.

\subsection{Description of the algorithm}\label{sect:algo}

We fix a total number $n$ of particles, as well as $1\in\left\{1,\ldots,n-1\right\}$ the number of resampled particles at each iteration of the algorithm.
In the sequel, when we consider a random variable $X_{i}^{j}$, the subscript $i$ denotes the index in $\left\{1,\ldots,n\right\}$ of a particle, and the superscript $j$ denotes the iteration of the algorithm.

In the algorithm below and in the following, we use classical notations for $k$-th order statistics. For $Y=(Y_1,\ldots,Y_n)$ an ensemble of independent and identically distributed (i.i.d.) real valued random variables with continuous cumulative distributions function, there exists  almost surely a unique (random) permutation $\sigma$ of $\left\{1,\ldots,n\right\}$ such that $Y_{\sigma(1)}<\ldots<Y_{\sigma(n)}$. For any $k \in \{1, \ldots,n\}$, we then use the classical notation $Y_{(k)}=Y_{\sigma(k)}$ to denote the $k$-th order statistics of the sample $Y$.

\medskip
 For any $x\in[0,a]$, we define the Adaptive Multilevel Splitting algorithm as follows (in order to approximate $p=\PP(X \ge a)$ one should take $x=0$, but we consider the general case $x\in[0,a]$ for theoretical purposes).

\begin{algo}[Adaptive Multilevel Splitting]\label{algo:AMS}
~

\noindent
{\bf Initialization:}
Define $Z^{0}=x$.
Sample $n$ i.i.d. realizations $X_{1}^{0},\ldots,X_{n}^{0}$, with the law $\mathcal{L}(X | X>x)$.

Define $Z^{1}=X_{(k)}^{0}$, the $k$-th order statistics of the sample $X^{0}=(X_{1}^{0},\ldots,X_{n}^{0})$, and $\sigma^1$ the (a.s.) unique associated permutation: $X_{\sigma^1(1)}^{0}<\ldots<X_{\sigma^1(n)}^{0}$. 

Set $j=1$.

\noindent
{\bf Iterations (on $j\geq 1$):}
While $Z^{j} <  a$:

Conditionally on $Z^{j}$, sample $k$ new independent \footnote{A precise mathematical statement is as follows. Let $(U_\ell^j)_{1\leq \ell\leq k,j\in\N^*}$ be i.i.d. random variables, uniformly distributed on $(0,1)$ and independent from all the other random variables. Then set $\chi_\ell^j=F(.;Z^j)^{-1}(U_\ell^j)$, where $F(.;x)^{-1}$ is the inverse distribution function associated with the (conditional) probability distribution $\mathcal{L}(X|X>x)$, see \ref{eq:fxy_cdf}. We slightly abuse notation by using $\mathcal{L}(X|X>Z^j)$ rather than $\mathcal{L}(X | X>z) |_{z=Z^{j}}$.} random variables $(\chi_1^j,\ldots,\chi_k^j)$, according to the law $\mathcal{L}(X | X>Z^{j})$. 


Set
$$
X_{i}^{j}=\begin{cases}\chi_{(\sigma^j)^{-1}(i)}^{j} \quad \text{if } (\sigma^j)^{-1}(i)\leq k\\ X_{i}^{j-1} \quad \text{if } (\sigma^j)^{-1}(i)>k. \end{cases}
$$

In other words, the particle with index $i$ is killed and resampled according to the law $\mathcal{L}(X | X>Z^{j})$ if $X_{i}^{j-1}\leq Z^{j}$, and remains unchanged if $X_{i}^{j-1}> Z^{j}$. Notice that the condition $(\sigma^j)^{-1}(i)\leq k$ is equivalent to $i\in\left\{\sigma^{j}(1),\ldots,\sigma^{j}(k)\right\}$.


Define $Z^{j+1}=X_{(k)}^{j}$, the $k$-th order statistics of the sample $X^{j}=(X_{1}^{j},\ldots,X_{n}^{j})$, and $\sigma^{j+1}$ the (a.s.) unique\footnote{The uniqueness of the permutation $\sigma^{j+1}$ is justified by Proposition \ref{propo:why_expo}.} associated permutation: $X_{\sigma^{j+1}(1)}^{j}<\ldots<X_{\sigma^{j+1}(n)}^{j}$. 

Finally increment $j\leftarrow j+1$.

\noindent
{\bf End of the algorithm:}
Define $J^{n,k}(x)=j-1$ as the (random) number of iterations. Notice that $J^{n,k}(x)$ is such that $Z^{J^{n,k}(x)} < a$ and $Z^{J^{n,k}(x)+1} \ge a$.

The estimator of the probability $p^x$ is defined by
\begin{equation}\label{eq:estimator}
\hat{p}^{n,k}(x)=C^{n,k}(x)\left(1-\frac{k}{n}\right)^{J^{n,k}(x)},
\end{equation}
with
\begin{equation}\label{eq:corrector}
C^{n,k}(x)=\frac{1}{n}{\rm Card}\left\{i ;\, X_{i}^{J^{n,k}(x)} \ge a\right\}.
\end{equation}
\end{algo}

Notice that $C^{n,1}(x)=1$. More generally, $C^{n,k}(x)\in\left\{\frac{n-k+1}{n},\ldots,\frac{n-k+k}{n}\right\}$.

Since we are interested in the algorithm starting at $x=0$, we introduce the notation
\begin{equation}\label{def0}
\hat{p}^{n,k}=\hat{p}^{n,k}(0).
\end{equation}

We finally stress that the computations of the sampled random variables $(X_i^0)_{1\leq i\leq n}$ for the initialization and of the $(\chi_{i}^{j})_{1\leq i\leq k}$ for each iteration $j$ can be made in parallel.

\subsection{Sampling from the conditional distributions $\mathcal{L}(X | X>z)$}\label{sect:samp}

At each iteration of the algorithm, we need to sample $k$ random variables according to conditional distributions $\mathcal{L}(X | X>z)$, with $z$ taking values in the sequence  $(Z^{j})_{0\leq j\leq J^{n,k}(x)}$.  As explained above, we develop our theoretical analysis of the properties of the algorithm (bias, variance and computational cost) in the idealized situation where it is possible to sample according to these conditional distributions $\mathcal{L}(X | X>z)$ for any $z\in[0,a]$. From a practical point of view, this assumption is generally unrealistic. One possible situation where it is realistic is the dynamic setting presented in the present section (in contrast with the static setting).

\subsubsection{The static setting and the exponential case}
\label{sec:static_expo}

In a general framework, there is no simple way to sample the distributions $\mathcal{L}(X | X>z)$. In practice, this can be done thanks to a Metropolis-Hastings procedure, see for example~\cite{guyader-hengartner-matzner-lober-11}. Of course, this introduces a bias and correlations between the particles at each iteration (compared with the idealized algorithm studied in this paper). This bias and these correlations asymptotically vanish when the number of iterations in the Metropolis-Hastings procedure goes to infinity. The error analysis associated to this procedure is out of the scope of this paper.

There is a simple example where it is actually possible to sample the distributions $\mathcal{L}(X | X>z)$, namely if $X$ is exponentially distributed. Indeed, if $X$ has exponential law $\mathcal{E}(1)$ with mean $1$, then the conditional distribution $\mathcal{L}(X | X>x)=\mathcal{L}(X+x)$ is a shifted exponential variable, for any $x>0$.
In the following, we will refer to this situation as {\em the exponential case}.
 Of course, this has no practical interest since in this case, $p=\PP(X \ge a) = \exp(-a)$ is analytically known. However, this particular case plays a crucial role in the analysis hereafter, since as will be precisely explained in Section~\ref{sect:wellposed}, the study of the general case can be reduced to the study of the exponential case after some change of variable. This trick was already used in the original papers~\cite{cerou-guyader-07a,cerou-guyader-07b}.




\subsubsection{Dynamic setting in dimension 1}
\label{sec:1ddyn}

In the one-dimensional dynamic setting, $X$ is defined as
\begin{equation*}
X=\sup_{0\leq t\leq \tau}Y_t
\end{equation*}
where $(Y_t)_{0\leq t\leq \tau}$  is a strongly Markovian time-homogeneous random process with values in $\R$, and $\tau$ is a stopping time. In this setting, the conditional distribution $\mathcal{L}(X |X>x)$ is easily sampled: it is the law of $\sup_{0 \leq t\leq \tau }Y_t^{x}$, where $Y_t^{x}$ denotes the stochastic process $(Y_t)_{ t \ge 0}$ which is such that $Y_0=x$.



Having in mind applications in molecular dynamics~\cite{cerou-guyader-lelievre-pommier-11}, a typical example is when $(Y_t)_{t\geq 0}$ satisfies a Stochastic Differential Equation
$$dY_t^x=f(Y_t^x)dt+\sqrt{2\beta^{-1}}dB_t, \quad Y_0^x=x,$$
with smooth drift coefficient $f$ and inverse temperature $\beta>0$. The stopping time is for example
$$\tau^x=\inf\left\{t\geq 0;Y_t^x<-\epsilon \text{ or } Y_t^x>1+\epsilon\right\},$$
for $x\in[0,1]$, and for some $\epsilon>0$ and one can then consider $X=\sup_{0\leq t\leq \tau^0}Y_t^0$. Let us consider the target level $a=1$. The AMS algorithm then yields an estimate of $\PP(X \ge 1)$ the probability that the stochastic process starting from $0$ reaches the level $1$ before the level $-\epsilon$. Such computations are crucial to compute transition rates and study the so-called reactive paths in the context of molecular dynamics, see~\cite{cerou-guyader-lelievre-pommier-11}.


Notice that in practice, a discretization scheme must be employed, which makes the sampling of the conditional probabilities more complicated. Another point of view on the AMS algorithm is then required in order to prove the unbiasedness of the estimator of the probability $p$, see~\cite{brehier-gazeau-goudenege-lelievre-rousset-2014}.




\begin{rem}
The exponential case can be obtained from a dynamic setting. Indeed, consider the following stochastic process: a particle starts at a given position $x$, moves on the real line with speed $p=+1$ on the random interval $[0,\tau]$ where $\tau\sim\mathcal{E}(1)$ is exponentially distributed, and with speed $p=-1$ on the interval $(\tau,+\infty)$. More precisely,
\begin{equation*}
p^x_t=\begin{cases} +1 \text{ for } 0\leq t\leq \tau\\ -1 \text{ for } t>\tau\end{cases}\quad 
q^x_t=\begin{cases} x+t \text{ for } 0\leq t\leq \tau\\ x+\tau-t \text{ for } t>\tau.\end{cases}
\end{equation*}
Notice that $(p^x_t,q^x_t)_{t \ge 0}$ is a Markov process such that $(q^x_t)_{t \ge 0}$ is continuous. Then for any initial condition~$x$ and any given threshold $a>x$ we have $\PP(\sup_{t\geq 0}q^x_t>a)=\PP(\tau>a-x)=\exp(x-a)$. In particular, $X=\sup_{t\geq 0}q^0_t=\sup_{0\leq t\leq \tau}q^0_t$ is an exponential random variable with parameter $1$.
\end{rem}

\subsubsection{Dynamic setting in higher dimension}
\label{sec:nddyn}

Let us consider again a strongly Markovian  time-homogeneous stochastic process $(Y_t)_{t\geq 0}$, but with values in $\mathbb{R}^d$ for $d\geq 2$. In this case, the levels need to be defined using a (continuous) function $\xi:\R^d\rightarrow \R$, sometimes called a reaction coordinate in the context of molecular dynamics.

Let us focus for simplicity on the case when $(Y_t)_{t\geq 0}$ is solution of the stochastic differential equation~(SDE):
\begin{equation}\label{SDEd}
dY_t^x=-\nabla V(Y_t^x)dt+\sqrt{2\beta^{-1}}dW_t, \quad Y_0^x=x,
\end{equation}
with smooth potential $V$, inverse temperature $\beta>0$ and $(W_t)_{t\geq 0}$ a $d$-dimensional Wiener process. 
Let us consider two disjoint closed subsets $A$ and $B$ of $\R^d$. Let us define the stopping time 
$$\tau^x=\min(\tau^x_A,\tau^x_B)$$
where
\begin{equation}\label{eq:tauAtauB}
\tau_A^x=\inf\left\{t\geq 0;Y_t^x \in A\right\} \text{ and }
\tau_B^x=\inf\left\{t\geq 0;Y_t^x \in B\right\}.
\end{equation}
Let us assume that the function $\xi$ is such that
$$A=\{x;\, \xi(x) \le 0\} \text{ and } B=\{x;\, \xi(x) \ge 1\}.$$
We then set, for a fixed initial condition $x_0 \in \R^d \setminus (A \cup B)$,
\begin{equation}\label{eq:dynamic}
X=\sup_{0\leq t\leq \tau}\xi(Y_{t}^{x_0}).
\end{equation}
Let us set $a=1$ as the target level.
In this case, the probability $p=\PP(X\geq 1)=\PP(\tau_B^{x_0}<\tau_A^{x_0})$ is the probability that the path starting from $x_0$ reaches~$B$ before~$A$. As explained above, this is a  problem of interest in molecular dynamics for example, to study reactive paths and compute transition  rates in high dimension, typically when~$A$ and~$B$ are metastable regions for $(Y_t)_{t \ge 0}$.



The problem to apply the AMS algorithm is again to sample according to the conditional distributions $\mathcal{L}(X | X>z)$. A natural idea is to use the following branching procedure in the resampling step at the $j$-th iteration: to build one of the new $k$ trajectories, one of the $(n-k)$ remaining trajectories is chosen at random, copied up to the first time it reaches the level $\{x; \, \xi(x)=Z^j\}$ and then completed independently from the past up to the stopping time $\tau$. The problem is that this yields in general a new trajectory which is correlated to the copied one through the initial condition on the level set $\{x; \, \xi(x)=Z^j\}$. Indeed, in general, given $x_1 \neq x_2$ such that $\xi(x_1)=\xi(x_2)$, the laws of $\sup_{0\leq t\leq \tau}\xi(Y_{t}^{x_1})$ and $\sup_{0\leq t\leq \tau}\xi(Y_{t}^{x_2})$ are not the same. As a consequence, it is unclear how to sample $\mathcal{L}(X | X>z)$, except if we would be able to build a function~$\xi$ such that the law of $\sup_{0\leq t\leq \tau}\xi(Y_{t}^{x})$ only depends on $\xi(x)$. This  is actually the case if $\xi$ is the so-called committor function associated to the dynamics~\eqref{SDEd} and the two sets $A$ and $B$.

\begin{defi}
Let $A$ and $B$ be two disjoint closed subsets in $\R^d$. The committor function $\xi$ associated with the SDE \eqref{SDEd} and the sets $A$ and $B$ is the unique solution of the following partial differential equation (PDE):
\begin{equation}\label{PDEcomm}
\begin{cases}
\begin{aligned}
-\nabla V \cdot \nabla\xi+\beta^{-1}\Delta \xi&=0 \quad \text{ in } \R^d \setminus (A\cup B),\\
\xi(x)&=0 \quad \text{ for } x\in A,\\
\xi(x)&=1 \quad \text{ for } x\in B.
\end{aligned}
\end{cases}
\end{equation}
\end{defi}

\begin{propo}\label{propo:comm}
Assume that $\xi$ is the committor function, solution of \eqref{PDEcomm}.
For any $x\in \R^d\setminus(A\cup B)$, we set $X^x=\sup_{0\leq t\leq \tau}\xi(Y_t^x)$. Notice that $X^x$ is a random variable with values in $[\xi(x),1]$.

We have $\PP(X^x > z)=\frac{\xi(x)}{z}$ for any $x\in\R^d\setminus(A\cup B)$ and $z \in [\xi(x), 1)$. In particular, if $\xi(x_1)=\xi(x_2)$, then we have the equality $\mathcal{L}(X^{x_1})=\mathcal{L}(X^{x_2})$.
Moreover, for any $x_1$, $x_2$ with $0<\xi(x_1)\leq \xi(x_2)<1$, $\mathcal{L}(X^{x_1}|X^{x_1}>\xi(x_2))=\mathcal{L}(X^{x_2})$.
\end{propo}

This previous Proposition (which is proven below) fully justifies the branching procedure described above to sample ${\mathcal L}(X^{x_0} | X^{x_0} > Z^j)$ at the $j$-th iteration of the algorithm: pick at random one of the $(n-k)$ remaining trajectories (say $(Y_t)_{t \ge 0}$), copy it up to the first time it reaches the level $\{x; \, \xi(x)=Z^j\}$ (let us denote $Y_\sigma$ the first hitting point of this level) and then complete the trajectory, independently from the past. By the strong Markov property, this yields a new $X$ sampled according to $\mathcal{L}(X^{Y_\sigma})$ which is indeed  $\mathcal{L}(X^{x_0} | X^{x_0} > Z^j)$, since $\xi(Y_\sigma)=Z^j$.

\begin{rem}
As already mentioned in the previous Section, in practice, the SDE~\eqref{SDEd} is discretized in time, say with a timestep $\Delta t$. Then, the branching procedure consists in copying the selected trajectory up to the first time index $n$ such that $\xi(Y_{n\Delta t}) > z$, and then to complete it independently from the past. This introduces a difference compared to the continuous in time situation considered above, since doing so we do not sample according to the conditional distribution $\mathcal{L}(X|X>z)$. To treat this situation, one needs to resort to other techniques to analyze the algorithm, see~\cite{brehier-gazeau-goudenege-lelievre-rousset-2014}. In particular, one can show that the algorithm still yields an unbiased estimator of $p$, using very different techniques of proof than those presented in this paper.  This approach is also useful to treat non-homogeneous in time Markov processes.
\end{rem}


At this stage, we can thus conclude that in the high-dimensional dynamic setting, if the committor function is known, the AMS algorithm can be practically implemented, and that it enters the framework of this paper. There are however two difficulties, that we will now discuss. 

First, the random variable $X=\sup_{0\leq t\leq \tau}\xi(Y^{x_0}_t)$, where $\xi$ is the committor function, does not satisfy Assumption~\ref{hyp:static}: we have $p=\PP(X\geq 1)=\PP(X=1)>0$ ($X$ takes values in $[0,1]$) and therefore, the cumulative distribution $t \mapsto \PP(X\le t)$ is not continuous at $t=1$. More precisely, from Proposition~\ref{propo:comm}, we have: $\forall t \in [0,1)$, $\PP(X\le t)=\left(1-\frac{\xi(x_0)}{t}\right)_+$ and $\forall t \ge 1$, $\PP(X\le t)=1$.
This is actually not a problem, as explained in the next Section in a general setting: the continuity of the cumulative distribution function is only required over $[0,1)$ (or more generally over $[0,a)$ in the general case when the target level is $a$).

The second difficulty is that knowing the committor function is actually a very strong assumption. Computing $\xi$ solution to~\eqref{PDEcomm} is actually impossible in practice since this is a high-dimensional PDE. Moreover, if $\xi$ was known, then we would actually know the small probability we want to estimate since $p=\PP(X \ge 1)=\xi(x_0)$. This is a consequence of the well-known probabilistic representation to solutions to~\eqref{PDEcomm}:
\begin{propo}
Recall the definitions~\eqref{eq:tauAtauB} of the stopping times $\tau_A^x$ and $\tau_B^x$. Then, if $\xi$ is the committor function associated to the SDE \eqref{SDEd} and the sets $A$ and $B$, then, for any $x\in\R^d$
\begin{equation}\label{eq:comm}
\xi(x)=\PP(\tau_B^x<\tau_A^x).
\end{equation}
\end{propo}


Thus, this high-dimensional dynamic case with known committor function should also be considered as an idealized setting, which is only useful for theoretical purposes, in order to study the best performance we could expect for the AMS algorithm.


We end up this Section with a proof of Proposition \ref{propo:comm}.

\noindent
{\bf Proof of Proposition \ref{propo:comm}:}
Let us consider $\xi$ satisfying~\eqref{PDEcomm}, $Y_t^x$ solution to~\eqref{SDEd} and  $X^x=\sup_{0\leq t\leq \tau}\xi(Y_t^x)$.
For any given $z\in(0,1)$, and any $x\in \R^d$, let us introduce
$$\tau_{z}^{x}=\inf\left\{t\geq 0;\xi(Y_{t}^{x})\geq z\right\}.$$
One easily checks the identity
$$\PP(X^x \ge z)=\PP(X^x>z)=\PP(\tau_{z}^{x}<\tau_{A}^{x}).$$

By continuity of $\xi$ and of the trajectories of the stochastic process $(Y_t)_{t\geq 0}$, and by the strong Markov property at the stopping time $\tau_z^x$, we get for any $x$ and any $z \in (\xi(x),1)$
\begin{align*}
\xi(x)&=\PP(\tau_{B}^{x}<\tau_{A}^{x})=\E\bigl[\mathds{1}_{\tau_{z}^{x}<\tau_{A}^{x}}\mathds{1}_{\tau_{B}^{x}<\tau_{A}^{x}}\bigr]\\
&=\E\bigl[\mathds{1}_{\tau_{z}^{x}<\tau_{A}^{x}}\mathds{1}_{\tau_{B}^{Y_{\tau_z^x}^x}<\tau_{A}^{Y_{\tau_z^x}^{x}}}\bigr]\\
&=\E\bigl[\mathds{1}_{\tau_{z}^{x}<\tau_{A}^{x}}\xi(Y_{\tau_z^x}^{x})\bigr]\\
&=\PP(\tau_{z}^{x}<\tau_{A}^{x})z.
\end{align*}


This identity proves the first claim of the Proposition. Moreover, since the law of $Y^x$ depends on $x$ only through $\xi(x)$, we also more generally get $\mathcal{L}(X^{x_1}|X^{x_1}>\xi(x_2))=\mathcal{L}(X^{x_2})$ for any $x_1$, $x_2$ with $0<\xi(x_1)\leq \xi(x_2)<1$,.
\qed

\subsection{About Assumption~\ref{hyp:static}}
\label{sec:hyp}

In this section, we show that Assumption~\ref{hyp:static} is actually too stringent. If one assumes the following
\begin{hyp}\label{hyp:static'}
$X$ is a real-valued positive random variable such that $t \in [0,a) \mapsto \PP(X \le t)$ is continuous,
\end{hyp}
\noindent
then the Algorithm~\ref{algo:AMS} is well defined, and  all the results presented below hold. In particular, the estimator $\hat{p}_{n,k}$ is an unbiased estimator of $$p=\PP(X \ge a).$$

We notice that Assumption~\ref{hyp:static'} is indeed more natural than Assumption~\ref{hyp:static} since the AMS algorithm only applies a resampling procedure with conditional distributions $\mathcal{L}(X|X>Z^j)$ to realizations such that $Z^j \in [0,a)$: this is why the continuity of the cumulative distribution function $t \mapsto \PP(X \le t)$ is actually only required over $[0,a)$.

The argument to show that one can recover the setting of Assumption~\ref{hyp:static} assuming only Assumption~\ref{hyp:static'} is the following coupling argument. In Lemma~\ref{lemma:pareto} below, it is proven that there exists a random variable $\tilde{X}$ such that:
\begin{itemize}
\item $\tilde{X}$ satisfies Assumption \ref{hyp:static};
\item for any $z\in [0,a)$, $\tilde{X} \le z$ is equivalent to  $X \le z$ and ${\mathcal L}(\tilde{X} | \tilde{X} >z)= {\mathcal L}(X | X >z)$;
\item $\tilde{X}\geq a$ is equivalent to $X \ge a$ and thus, in particular, $\PP(\tilde{X}\geq a)=\PP(X \ge a)=p$, the probability to be estimated.
\end{itemize}
The last two properties show that running the AMS algorithm on $\tilde{X}$ is equivalent to running the AMS algorithm on $X$: the iterations, the stopping criterion and the estimator are the same.  
The theory developed in this paper (unbiased estimator, analysis of the cost and of the computational cost) is then applied to the algorithm applied to the auxiliary random variable $\tilde{X}$ instead of $X$, which is completely equivalent to the algorithm applied to $X$.


In all the following, for simplicity, we will always assume that Assumption~\ref{hyp:static} holds, keeping in mind that it can be relaxed to Assumption~\ref{hyp:static'}. Thus, inequalities which involve the random variable $X$ can be changed from large to strict without modifying the associated events (almost surely).

We end this Section with a Lemma which defines the random variable $\tilde{X}$ as a function of $X$.
\begin{lemme}\label{lemma:pareto}
Let $X$ be a random variable satisfying Assumption~\ref{hyp:static'}, and let us define
$$\tilde{X}=X\mathds{1}_{X<a}+\frac{a}{U}\mathds{1}_{X\ge a},$$
where $U$ is a random variable independent of $X$ and uniformly distributed on $(0,1)$.

Then, (i) $\tilde{X}$ satisfies Assumption \ref{hyp:static}, (ii) for any $z\in [0,a)$, $\tilde{X} \le z$ is equivalent to  $X \le z$, and the two laws ${\mathcal L}(\tilde{X} | \tilde{X} >z)$ and ${\mathcal L}(X | X >z)$ coincide on $(z,a)$ and (iii) $\tilde{X}\geq a$ is equivalent to $X \ge a$.
\end{lemme}


\noindent
{\bf Proof of Lemma \ref{lemma:pareto}:}
Since $a/U>a$, it is easy to check the items (ii) and (iii). Let us now consider the cumulative distribution of $\tilde{X}$.

For $t<a$, $\PP(\tilde{X} \le t) = \PP(X \le t)$ and thus $t \mapsto \PP(\tilde{X} \le t)$ is continuous for $t \in [0,a)$ by Assumption~\ref{hyp:static'}.

For $t>a$, $\PP(\tilde{X} \le t) = \PP(X < a) + \PP( U \ge a/t, \, X \ge a)= \PP(X<a) + \PP(X \ge a) \left(1-\frac{a}{t} \right) = 1 - \frac{a}{t} \PP(X \ge a)$ and thus $t \mapsto \PP(\tilde{X} \le t)$ is continuous for $t \in (a,+\infty)$.

Finally, with these expressions one easily checks left and right continuity of $t\mapsto\PP(\tilde{X} \le t)$ at $a$.

This concludes the proof of the fact that $\tilde{X}$ satisfies Assumption \ref{hyp:static}, and thus the proof of Lemma~\ref{lemma:pareto}.
\qed


\section{Reduction to the exponential case and well-posedness of the algorithm}
\label{sec:expo_wp}

The aim of this Section (see Section~\ref{sect:wellposed}) is to prove the well-posedness of the algorithm, namely the fact that $J^{n,k}(x)$ is almost surely finite, when the probability $p=P(0)$ to estimate is positive. The argument is based on the fact that the general case is related to the exponential case through a change of variable, see Section~\ref{sec:expo} (this will be instrumental in the rest of the paper). Section~\ref{sect:notations} first gives a few notation that will be useful below.

\subsection{Notation}\label{sect:notations}

\subsubsection{General notation}



We will use the following set of notations, associated to the random variable $X$ satisfying Assumption~\ref{hyp:static}.

We denote by $F$  the cumulative distribution function of the random variable $X$: $F(t)=\PP(X\leq t)$ for any $t\in \R$, and $F(0)=0$. From Assumption~\ref{hyp:static}, the function $F$ is continuous. Notice that it ensures that if $Y$ is an independent copy of $X$, then $\PP(X=Y)=0$, and thus, in the algorithm, there is only at most one sample at a given level.


We recall that our aim is to estimate the probability
\begin{equation}
p=\PP(X \ge a)=\PP(X>a),
\end{equation}
given a threshold $a>0$. More generally, we define for any $x\in[0,a]$
\begin{equation}
P(x)=\PP(X>a | X>x),
\end{equation}
so that we have $p=P(0)$. Notice that $P(a)=1$.

For any $x\in[0,a]$, $\mathcal{L}(X | X>x)$ admits a cumulative distribution function $F(.;x)$, which satisfies: for any $y \in \R$,
\begin{equation}\label{eq:fxy_cdf}
\begin{gathered}
F(y;x)=\frac{F(y)-F(x)}{1-F(x)}\mathds{1}_{y\geq x}.
\end{gathered}
\end{equation}
While the variable $x$ always denotes the parameter in the conditional distributions $\mathcal{L}(X | X>x)$, we use the variable $y$ as a dummy variable in the associated densities and cumulative distribution functions.

By Assumption \ref{hyp:static}, $\mathds{1}_{y\geq x}$ in the definition above can be replaced with $\mathds{1}_{y>x}$. Notice that $F(y;0)=F(y)$. Moreover, with these notations, we have
\begin{equation}\label{eq:P}
P(x)=1-F(a;x).
\end{equation}

An important tool in the following is the family of functions: for any $x\in[0,a]$ and any $y\in \R$
\begin{equation}\label{def:Lambda}
\begin{gathered}
\Lambda(y;x)=-\log(1-F(y;x))\in[0,+\infty]\\
\Lambda(y)=\Lambda(y;0)=-\log(1-F(y)).
\end{gathered}
\end{equation}

We remark the following identity: for $0\leq x\leq y\leq a$,
$$\Lambda(y;x)=\log(1-F(x))-\log(1-F(y))=\Lambda(y)-\Lambda(x).$$

\subsubsection{Specific notation when $X$ admits a density $f$}
\label{sec:density}

In some places, we will assume that $X$ admits a density $f$ with respect to the Lebesgue measure (which indeed implies Assumption~\ref{hyp:static}). This assumption is in particular satisfied in the exponential case (namely when $X$ is exponentially distributed), which we will consider in several arguments below to study the bias and the computational cost. 

If $X$ admits a density $f$,  the law $\mathcal{L}(X | X>x)$ of $X$ conditionally on $\left\{X>x\right\}$ also admits a density $f(.;x)$, which satisfies: for any $y\geq 0$,
\begin{equation}\label{eq:fxy_density}
\begin{gathered}
f(y;x)=\frac{f(y)}{1-F(x)}\mathds{1}_{y\geq x}.
\end{gathered}
\end{equation}
 Notice that $f(y)=f(y;0)$. 

We finally introduce some notations about order statistics of samples of i.i.d. real random variables. If $X$ admits a density $f$, then the $k$-th order statistics of an i.i.d. sample $(X_1, \dots , X_n)$ (distributed according to the law of $X$) admits a density $f_{n,k}$ which satisfies: for any $y\in \R$,
$$f_{n,k}(y)=k\binom{n}{k}F(y)^{k-1}f(y)\bigl(1-F(y)\bigr)^{n-k}.$$
The associated cumulative distribution function is $F_{n,k}(y)=\int_{-\infty}^{y}f_{n,k}(z)dz$.
Likewise, we introduce notations for the density and the cumulative distribution function of order statistics for the law $\mathcal{L}(X | X>x)$: when $0\leq x\leq y\leq a$ we set
\begin{equation}\label{statordre}
\begin{gathered}
f_{n,k}(y;x)=k\binom{n}{k}F(y;x)^{k-1}f(y;x)\bigl(1-F(y;x)\bigr)^{n-k},\\
F_{n,k}(y;x)=\int_{x}^{y}f_{n,k}(z;x)dz.
\end{gathered}
\end{equation}

\subsection{Reduction to the exponential case}
\label{sec:expo}

One of the key tool in the following is the reduction of the general case to the exponential case, thanks to the use of the function $\Lambda$, defined by \eqref{def:Lambda}. We recall that the exponential case refers to the case when $X$ is distributed according to the exponential law with parameter 1, see Section~\ref{sec:static_expo}. The basic remark is the following. Since by Assumption \ref{hyp:static} the cumulative distribution function $F(.;x)$ is continuous, we have the following classical result, see for instance Proposition $2.2$ in \cite{asmussen-glynn}.
\begin{lemme}\label{lemma:Lambda_expo}
If $Y\sim \mathcal{L}(X | X>x)$, then $F(Y;x)$ is uniformly distributed on $(0,1)$, and thus $\Lambda(Y;x)$ has an exponential law with parameter $1$.
\end{lemme}

Let us first state a result on the algorithm without any stopping criterion.

\begin{propo}\label{propo:why_expo}
Let us consider the sequence of random variables $((X_i^j)_{1 \le i \le n},Z^j)_{j \ge 0}$ generated by the AMS algorithm~\ref{algo:AMS} {\em  without any stopping criterium}.
Set $Y_{i}^{j}=\Lambda(X_{i}^{j})$ and $S^{j}=\Lambda(Z^j)$. Then we have the following properties.
\begin{itemize}
\item[(i)] For any $j\geq 0$, $(Y_{i}^{j}-S^j)_{1\leq i\leq n}$ is a family of i.i.d. exponentially distributed random variables, with parameter $1$.
\item[(ii)] For any $j\geq 1$, $(Y_{i}^{j}-S^j)_{1\leq i\leq n}$ is independent of $(S^l-S^{l-1})_{1\leq l\leq j}$.
\item[(iii)] The sequence $(S^{j}-S^{j-1})_{j\geq 1}$ is i.i.d..
\end{itemize}
As a consequence, in law, the sequence $((Y_{i}^{j})_{1\leq i\leq n}, S^j)_{j\geq 0}$ is equal to the sequence of random variables $((X_i^j)_{1 \le i \le n},Z^j)_{j \ge 0}$ obtained by the realization of the AMS algorithm {\em without any stopping criterion} in the exponential case, with initial condition $Z^0=\Lambda(x)$.
\end{propo}

\noindent
{\bf Proof of Proposition \ref{propo:why_expo}:}
The last assertion is a direct consequence of the three former items and of the fact that in the exponential case the function $\Lambda$ is the identity mapping.

Item (iii) is a consequence of items (i) and (ii). It remains to prove jointly those two items, which is done by induction on $j$. When $j=0$, the result follows from the way the algorithm is initialized: for each $1\leq i\leq n$, $Y_{i}^{0}=\Lambda(X_{i}^{0})$ is exponentially distributed (see Lemma \ref{lemma:Lambda_expo}), and the independence property is clear. Assuming that the properties (i) and (ii) are satisfied for all $k\leq j$, it is sufficient to prove that:
\begin{enumerate}[(i)]
\item $(Y_{i}^{j+1}-S^{j+1})_{1\leq i\leq n}$ are i.i.d. and exponentially distributed with mean $1$;
\item $(Y_{i}^{j+1}-S^{j+1})_{1\leq i\leq n}$ is independent of $(S^l-S^{l-1})_{1\leq l\leq j+1}$.
\end{enumerate}
For any positive real numbers $y_1,\ldots,y_n$, and any $s^1,\ldots,s^{j+1}$, this is equivalent to proving that
\begin{equation}\label{eq:induction}
\begin{aligned}
A:&=\PP\left(Y_{1}^{j+1}-S^{j+1}>y_1,\ldots,Y_{n}^{j+1}-S^{j+1}>y_n,S^{1}-S^{0}>s^1,\ldots,S^{j+1}-S^{j}>s^{j+1}\right)\\
&=\exp\left(-(y_1+\dotsb+y_n)\right)\PP(S^{1}-S^{0}>s^1,\ldots,S^{j+1}-S^{j}>s^{j+1}).
\end{aligned}
\end{equation}

We want to decompose the probability with respect to the value of $\bigl(\sigma^{j+1}(i))_{1\leq i\leq k}$. We recall that almost surely we have
$$Y_{\sigma^{j+1}(1)}^{j}<\ldots<Y_{\sigma^{j+1}(k)}^{j}=S^{j+1}<Y_{\sigma^{j+1}(k+1)}^{j}<\ldots<Y_{\sigma^{j+1}(n)}^{j}.$$
In fact, in order to preserve symmetry inside the groups of resampled and not-resampled particles, we decompose over the possible values for the random set $\left\{\sigma^{j+1}(1) \ldots \sigma^{j+1}(k)  \right\}$. We thus compute a sum over all partitions
$$\left\{1,\ldots n\right\}=I_{-}\sqcup I_{+},$$
such that $\text{Card}\bigl(I_{-}\bigr)=k$.


\begin{align*}
A&=\PP\Bigl(Y_{1}^{j+1}-S^{j+1}>y_1,\ldots,Y_{n}^{j+1}-S^{j+1}>y_n,S^{1}-S^{0}>s^1,\ldots,S^{j+1}-S^{j}>s^{j+1}\Bigr)\\
&=\sum_{\overset{I_{-}\subset\left\{1,\ldots,n\right\}}{\text{Card}(I_{-})=k}}\PP\Bigl( Y_{1}^{j+1}-S^{j+1}>y_1,\ldots,Y_{n}^{j+1}-S^{j+1}>y_n,\\
& \hspace{2cm} S^{1}-S^{0}>s^1,\ldots,S^{j+1}-S^{j}>s^{j+1},\left\{\sigma^{j+1}(1), \ldots , \sigma^{j+1}(k) \right\}=I_{-}\Bigr)\\
&=\sum_{\overset{I_{-}\subset\left\{1,\ldots,n\right\}}{\text{Card}(I_{-})=k}}\PP\Bigl( \left\{Y_{i}^{j+1}-S^{j+1}>y_{i}; i\in I_{+}\right\},S^{1}-S^{0}>s^1,\ldots,S^{j+1}-S^{j}>s^{j+1}\Bigr)\\
&\hspace{2cm}\prod_{i\in I_{-}}\PP\Bigl(\Lambda(\chi^{j+1})- \Lambda(Z^{j+1}) >y_i \Bigr).\\
\end{align*}
In the last line, we used the fact that by construction of the algorithm, on the event we consider, namely $\left\{ \sigma^{j+1}(1), \ldots , \sigma^{j+1}(k) \right\}=I_{-}$, we have the equality of random variables $Y^{j+1}_i = \Lambda\left(\chi_{(\sigma^{j+1})^{-1}(i)}^{j+1}\right), i \in I_-$. The latter are i.i.d. and independent of all the other random variables used at this stage of the algorithm. Moreover:
$$\prod_{i\in I_{-}}\PP\Bigl(\Lambda(\chi^{j+1})- \Lambda( Z^{j+1} ) > y_i\Bigr)=\prod_{i\in I_{-}}\exp\bigl(-y_i\bigr).$$

Let us now introduce a notation: for $I_{-} \subset \{ 1, \cdots n \}$ we set
\[
M^{j}_{I_-}:=\max\left\{Y_{i}^{j}-S^j; i\in I_{-}\right\}.
\]
Remark that on the event $\left\{ \sigma^{j+1}(1), \ldots , \sigma^{j+1}(k) \right\}=I_{-}$, almost surely, we have $Y^{j}_i = Y^{j+1}_i$ and $S^{j+1}-S^j = M^{j}_{I_-}$. Using the independence properties of the $Y_{i}^{j}-S^j, 1\leq i \leq n $ (and thus of $M^{j}_{I_-}$) from the induction hypothesis $(ii)$, we obtain
\begin{align*}
\PP\Bigl(&\left\{ Y^{j+1}-S^{j+1}>y_{i}; i\in I_{+}\right\},S^{1}-S^{0}>s^1,\ldots,S^{j+1}-S^{j}>s^{j+1}\Bigr)\\
&=\PP\Bigl(\left\{Y_{i}^{j}-S^{j}-M^{j}_{I_-}>y_{i}; i \in I_{+}\right\},S^{1}-S^{0}>s^1,\ldots, M^{j}_{I_-}>s^{j+1}\Bigr)\\
&=\PP\Bigl(\left\{Y_{i}^{j}-S^{j}-M^{j}_{I_-}>y_{i}; i \in I_{+}\right\},  M^{j}_{I_-} >s^{j+1}\Bigr)\PP\Bigl(S^{1}-S^{0}>s^1,\ldots,S^{j}-S^{j-1}>s^{j}\Bigr).
\end{align*}
We can then integrate the $Y_{i}^{j}-S^j,  i \in I_+$ using the induction hypothesis $(i)$: 
\begin{align*}
\PP\Bigl(&\left\{Y_{i}^{j}-S^j-M^{j}_{I_-}>y_{i} ; i\in I_{+} \right\},M^{j}_{I_-}>s^{j+1}\Bigr)\\
&=\E\Bigl[\prod_{i\in I_{+}}\exp\Bigl(-y_{i}-M^j_{I_-}\Bigr)\mathds{1}_{M^j_{I_-}>s^{j+1}}\Bigr]\\
&=\prod_{i\in I_{+}}\exp\bigl(-y_i\bigr)\E\Bigl[\exp\Bigl(- (n-k)M^j_{I_-} \Bigr)\mathds{1}_{M^j_{I_-}>s^{j+1}}\Bigr]
\end{align*}

The proof is now complete since:
\begin{align*}
A&=\PP\Bigl(Y_{1}^{j+1}-S^{j+1}>y_1,\ldots,Y_{n}^{j+1}-S^{j+1}>y_n,S^{1}-S^{0}>s^1,\ldots,S^{j+1}-S^{j}>s^{j+1}\Bigr)\\
&=\sum_{\overset{I_{-}\subset\left\{1,\ldots,n\right\}}{\text{Card}(I_{-})=k}}\Bigl(\prod_{i\in I_{-}}\exp\bigl(-y_i\bigr)\Bigr)\Bigl(\prod_{i\notin I_{-}}\exp\bigl(-y_i\bigr)\Bigr)\E\Bigl[\exp\Bigl(-(n-k)M^j_{I_-}\Bigr)\mathds{1}_{M^j_{I_-}>s^{j+1}}\Bigr]\\
&\hspace{3cm}\PP\Bigl(S^{1}-S^{0}>s^1,\ldots,S^{j}-S^{j-1}>s^{j}\Bigr)\\
&=\exp\Bigl(\sum_{i=1}^{n}y_i\Bigr)\PP\Bigl(S^{1}-S^{0}>s^1,\ldots,S^{j}-S^{j-1}>s^{j}\Bigr)\sum_{\overset{I_{-}\subset\left\{1,\ldots,n\right\}}{\text{Card}(I_{-})=k}}\E\Bigl[\exp\Bigl(-(n-k) M^j_{I_-} \Bigr)\mathds{1}_{M^j_{I_-}>s^{j+1}}\Bigr],
\end{align*}
which proves \eqref{eq:induction}.

In particular, taking $y_1=\ldots=y_n=s^1=\ldots =s^j=0$, we see that
$$\sum_{\overset{I_{-}\subset\left\{1,\ldots,n\right\}}{\text{Card}(I_{-})=k}}\E\Bigl[\exp\Bigl( - (n-k) M^j_{I_-}\Bigr)\mathds{1}_{M^j_{I_-} > s^{j+1}}\Bigr]=\PP(S^{j+1}-S^{j}>s^{j+1}).$$

This concludes the proof of Proposition \ref{propo:why_expo}.

\qed

The next Lemma shows that the AMS algorithm applied to $X$ with target level $a$ and the AMS algorithm applied to $\Lambda(X)$ with target level $\Lambda(a)$ stop at the same iteration.
\begin{lemme}\label{lemme:as_equal}
Let us consider the sequence of random variables $((X_i^j)_{1 \le i \le n},Z^j)_{j \ge 0}$  generated by the AMS algorithm~\ref{algo:AMS} {\em without any stopping criterium},  and set $Y_{i}^{j}=\Lambda(X_{i}^{j})$ and $S^{j}=\Lambda(Z^j)$. For any $\alpha > 0$,
almost surely, $\left\{S^j\geq\Lambda(\alpha)\right\}=\left\{Z^j\geq\alpha\right\}$.
\end{lemme}
\noindent
{\bf Proof of Lemma \ref{lemme:as_equal}:} : First, $\Lambda$ is non-decreasing so that $\left\{Z^j\geq\alpha\right\}\subset\left\{S^j\geq\Lambda(\alpha)\right\}$.
Moreover, one easily checks that $\left\{S^j\geq\Lambda(\alpha)\right\}\cap\left\{Z^j<\alpha\right\}\subset \left\{S^j=\Lambda(\alpha)\right\}$. From Proposition \ref{propo:why_expo} we know that $S^j=S^0+\sum_{\ell=1}^{j}\left(S^{\ell}-S^{\ell-1}\right)$ admits a density with respect to the Lebesgue measure, since the $S^{\ell}-S^{\ell-1}$ have the density of the $k$-th order statistics of independent exponentially distributed random variables with parameter~$1$. Therefore $\PP\left(\left\{S^j>\Lambda(\alpha)\right\}\neq\left\{Z^j>\alpha\right\}\right)=0$.
\qed

A direct corollary of Proposition~\ref{propo:why_expo} and Lemma~\ref{lemme:as_equal} is that the original problem reduces to the exponential case. 
\begin{cor}\label{cor:expo}
Consider the sequences of random variables $(X_i^j)_{1 \le i \le n, 0 \le j \le J^{n,k}(x)}$ and~$(Z^j)_{0 \le j \le J^{n,k}(x)+1}$ generated by the AMS algorithm~\ref{algo:AMS}. Set $Y_{i}^{j}=\Lambda(X_{i}^{j})$ and $S^{j}=\Lambda(Z^j)$.

Then, the sequences $(Y_{i}^{j})_{1\leq i\leq n, 0 \le j \le J^{n,k}(x)}$ and $(S^j)_{0 \le j \le J^{n,k}(x)+1 }$ are equal in law to the sequences of random variables $(X_i^j)_{1 \le i \le  J^{n,k}(\Lambda(x))}$ and $(Z^j)_{0 \le j \le J^{n,k}(\Lambda(x))+1}$ obtained by the realization of the AMS algorithm in the exponential case, with initial condition $Z^0=\Lambda(x)$ and target level $\Lambda(a)$.
\end{cor}

Finally, in the next Sections, we need the following result, which is a consequence of Proposition \ref{propo:why_expo}.
\begin{cor}\label{Cor:condit}
For any $j\geq 0$, conditionally on $Z^{j}$, the random variables $(X_{i}^{j})_{1\leq i\leq n}$ are i.i.d. with law $\mathcal{L}(X | X>Z^{j})$.
\end{cor}
\noindent
{\bf Proof:}
Thanks to Proposition \ref{propo:why_expo}, we see that $(\Lambda(X_{i}^{j})-\Lambda(Z^{j}))_{1\leq i\leq n}$ are i.i.d. and exponentially distributed with mean $1$. Since $\Lambda(X_{i}^{j})-\Lambda(Z^{j})=\Lambda(X_{i}^{j};Z^{j})$, we observe that for any $x_1,\ldots,x_n\in [Z^j,+\infty)$,
\begin{align*}
\PP(X_{1}^{j} > x_1&,\ldots,X_{n}^{j}> x_n | Z^j)=\PP(\Lambda(X_{1}^{j})>\Lambda(x_1),\ldots,\Lambda(X_{n}^{j})>\Lambda(x_n) | \Lambda(Z^j))\\
&=\PP(\Lambda(X_{1}^{j})-\Lambda(Z^j)>\Lambda(x_1)-\Lambda(Z^j),\ldots,\Lambda(X_{n}^{j})-\Lambda(Z^j)>\Lambda(x_n)-\Lambda(Z^j) | \Lambda(Z^j))\\
&=\exp(-(\Lambda(x_1)-\Lambda(Z^j))\ldots \exp(-(\Lambda(x_n)-\Lambda(Z^j))\\
&=(1-F(x_1;Z^j))\ldots (1-F(x_n;Z^j)).
\end{align*}
This concludes the proof.
\qed

\subsection{Well-posedness of the algorithm}\label{sect:wellposed}
To ensure that the algorithm giving an estimator of the probability is well-defined, namely that it gives a result after a finite number of steps, we prove in this Section that $J^{n,k}(x)$ is almost surely finite, when the probability $p=P(0)$ is positive. The proof relies on the reduction to the exponential case explained in the previous Section.
\begin{propo}\label{prop:finiteJ}
Suppose $p=P(0)>0$. Then for any $x\in[0,a]$, we have $P(x)=\PP(X>a | X>x)>0$, and for any integers $n$ and $k$ with $1\leq k<n$, the number of iterations in the AMS algorithm is almost surely finite: $J^{n,k}(x)<+\infty$ a.s.
\end{propo}

\noindent
{\bf Proof of Proposition \ref{prop:finiteJ}:}
To prove this result, we consider the AMS Algorithm~\ref{algo:AMS} without any stopping criterion (namely the condition $Z^j \geq a$). As a consequence, we define sequences of random variables with the iteration index $j\in \N$: we get $(X_{i}^{j})_{j\geq 0}$, for any $i \in \{1,\ldots, n\}$ and $(Z^{j})_{j\geq 0}$. Proposition \ref{prop:finiteJ} is then equivalent to the following statement: almost surely, $\left\{j\geq 0; Z^j \ge a\right\}\neq \emptyset$.
 
Thanks to Proposition \ref{propo:why_expo}, we write for any $j\geq 0$
$$S^j=S^0+\sum_{\ell=1}^{j}R^{\ell},$$
where $R^{\ell}=S^{\ell}-S^{\ell-1}$ are independent and identically distributed positive random variables, satisfying $\E R^{\ell}\in(0;+\infty)$. Indeed,
$$\E R^1\leq \E\max_{1\leq i\leq n}Y_{i}^{0}\leq \E \sum_{1\leq i\leq n}Y_{i}^{0}=n,$$
where we recall that $Y_{i}^{0}=\Lambda(X_{i}^{0})$ are independent and exponentially distributed with parameter $1$. To prove that $\E R^{1}>0$, we write
$$\E R^1\geq \E\min_{1\leq i\leq n}Y_{i}^{0}=1/n,$$
since it is easily checked that $\min_{1\leq i\leq n}Y_{i}^{0}$ has an exponential distribution, with mean $1/n$.

By the Strong Law of Large Numbers, when $j\rightarrow +\infty$, we have the almost sure convergence
$$\frac{S^j}{j}\rightarrow \E R^1,$$
which yields $S^j\rightarrow +\infty$, almost surely, when $j\rightarrow +\infty$. As a consequence, almost surely, there exists some $j\in \N$ such that $S^j \ge \Lambda(a)$. Using Lemma~\ref{lemme:as_equal}, this then implies that $Z^j \ge a$, and that $J^{n,k}(x)<+\infty$.
\qed

In the case $k=1$, following the ideas in the proof of  Proposition \ref{prop:finiteJ}, one can easily identify the law of the number of iterations (see~\cite{guyader-hengartner-matzner-lober-11} for a similar result).
\begin{propo}\label{propo:Poisson1}
The random variable $J^{n,1}(x)$ has a Poisson distribution with mean $-n\log(P(x))$.
\end{propo}

\noindent
{\bf Proof of Proposition \ref{propo:Poisson1}:} In the case $k=1$, $R^1=\min_{1\leq i\leq n}Y_{i}^{0}$ has an exponential distribution with mean $1/n$. We recall that $(R^{\ell}=S^{\ell}-S^{\ell-1})_{\ell\geq 0}$ is a sequence of independent and identically distributed random variables.

Let us introduce the Poisson process, with intensity $n$, associated with the sequence of independent and exponentially distributed increments $(R^{j})_{\ell\geq 0}$:
$$P_t=\sum_{\ell=0}^{+\infty}\mathds{1}_{S^\ell \leq t}.$$

Since $S^0=\Lambda(x)$, we identify that
$$J^{n,1}(x)=P_{\Lambda(a)-\Lambda(x)}=P_{-\log(P(x))}.$$
The result follows since for any $t\geq 0$ $P_t$ has a Poisson distribution with mean $nt$.
\qed

\section{The estimator $\hat{p}^{n,k}(x)$ is unbiased}\label{sect:unbiased}

Recall that $X$ satisfies Assumption \ref{hyp:static}.

Let us fix a total number of replicas $n$, as well as $k \in \{1, \ldots,n-1\}$ the number of killed and resampled replicas at each iteration. Given a level $a>0$, we recall that the estimator of the conditional probability $P(x)=\PP(X>a | X>x)$ for each value of $x\in[0,a]$ is $\hat{p}^{n,k}(x)$, defined by \eqref{eq:estimator}. We introduce the following notation:
\begin{equation}\label{eq:bias}
p^{n,k}(x)=\E[\hat{p}^{n,k}(x)].
\end{equation}

Recall that we are specifically interested in estimating the probability $p=P(0)$, and that the introduction of $P(x)$ for $x$ in the interval $[0,a]$ is a tool to prove that the estimator is unbiased. We write the result in its full generality, and then specify it to the estimation of $p$.

\begin{theo}\label{th:unbiased}
For any $k \in \{1,\ldots, n-1\}$, for any $a>0$, such that $p=\PP(X>a)=P(0)>0$, and any $x\in[0,a]$, $\hat{p}^{n,k}(x)$ is an unbiased estimator of the conditional probability $P(x)$:
\begin{equation}\label{eq:nobias}
\E[\hat{p}^{n,k}(x)]=P(x).
\end{equation}

In particular, when $x=0$, we have $\E[\hat{p}^{n,k}]=p$.
\end{theo}

From Section~\ref{sec:expo} (see Corollary~\ref{cor:expo}), it is sufficient to prove the result in the exponential case. Indeed, let us assume that~\eqref{eq:nobias} holds in the exponential case, and let us consider the general case of a random variable $X$ satisfying Assumption~\ref{hyp:static}, then we have
\begin{align*}
\E[\hat{p}^{n,k}(x)]
&= \PP(\Lambda(X) > \Lambda(a) | \Lambda(X) > \Lambda(x)) \\
&= \exp(-\Lambda(a) + \Lambda(x))\\
&= \PP(X > a | X > x)=P(x).
\end{align*}
The first equality is a consequence of Corollary~\ref{cor:expo} and Theorem~\ref{th:unbiased} in the exponential case. The third equality is a direct consequence of the definition~\eqref{def:Lambda} of $\Lambda$.

The aim of this section is thus to prove Theorem~\ref{th:unbiased} in the exponential case. In all the following, we denote by $f(x)=\exp(-x)$ the density of $X$, and we will use the notation introduced in Section~\ref{sec:density} above for the density of the $k$-th statistics. Actually, the proof given below is valid as soon as $X$ has a density $f$: we do not use the specific form of the density, and this specific form would not make the argument easier. 

The proof of this result is divided into two steps. First, we show that the function $x\mapsto p^{n,k}(x)$ is solution of a functional equation. Second, we show that the function $x\mapsto P(x)$ is its unique solution.

\subsection{Derivation of the functional equation satisfied by $p^{n,k}$}

\begin{propo}\label{propo:funct_p}
Let us assume that $X$ admits a density $f$. Assume $p=P(0)>0$. The function \linebreak $x\in[0,a]\mapsto p^{n,k}(x)$ is solution of the following functional equation (with unknown $q$): for any $0\leq x\leq a$
\begin{equation}\label{funct_p}
q(x)=\int_{x}^{a}\left(1-\frac{k}{n}\right)q(y)f_{n,k}(y;x)dy+\theta_{p}^{n,k}(x),
\end{equation}
with
\begin{equation}\label{eq:theta_p}
\theta_{p}^{n,k}(x)=\sum_{l=0}^{k-1}\frac{n-l}{n}\PP\left(S(x)_{(l)}^{n}\leq a\leq S(x)_{(l+1)}^{n}\right),
\end{equation}
where $(S(x)_{l}^{n})_{1\leq l\leq n}$ are independent and identically distributed with density $f(.;x)$ (see~\eqref{statordre}), while for $1\leq l \leq k-1$, $S(x)_{(l)}^{n}$ denotes the $l$-th order statistics of this $n$-sample. By convention, we set $S(x)_{(0)}^{n}=x$.
\end{propo}
\noindent
{\bf Proof:}
The key idea is to decompose the expectation $\E[\hat{p}^{n,k}]$ according to the (random but almost surely finite) value of the number $J^{n,k}(x)$ of iterations. The function $\theta_{p}^{n,k}$ appears as the result of the algorithm when $J^{n,k}(x)=0$, while the integral formulation corresponds to the case $J^{n,k}(x)>0$. In the latter case, we then condition on the value of the first level $Z^{1}=X_{(k)}^{0}$ and use Corollary \ref{Cor:condit}.

More precisely, we have
$$p^{n,k}(x)=\E\left[\hat{p}^{n,k}(x)\right]=\E\left[\hat{p}^{n,k}(x)\mathds{1}_{J^{n,k}(x)=0}\right]+\E\left[\hat{p}^{n,k}(x)\mathds{1}_{J^{n,k}(x)>0}\right].$$

First, we have from \eqref{eq:estimator} and \eqref{eq:corrector}
\begin{align*}
\E\left[\hat{p}^{n,k}(x)\mathds{1}_{J^{n,k}(x)=0}\right]&=\E\left[C^{n,k}(x)\mathds{1}_{J^{n,k}(x)=0}\right]=\E\left[C^{n,k}(x)\mathds{1}_{X_{(k)}^{0}>a}\right]\\
&=\sum_{l=0}^{k-1}\frac{n-l}{n}\E\left[\mathds{1}_{X_{(l)}^{0}\leq a <X_{(l+1)}^{0}}\right]=\theta_{p}^{n,k}(x).
\end{align*}

Second, using conditional expectation with respect to $Z^1$,
\begin{align*}
\E\left[\hat{p}^{n,k}(x)\mathds{1}_{J^{n,k}(x)>0}\right]&=\E\left[\E\left[\left(1-\frac{k}{n}\right)\left(1-\frac{k}{n}\right)^{J^{n,k}(x)-1}C^{n,k}(x) | Z^1\right]\mathds{1}_{Z^1\leq a}\right]\\
&=\E\left[\left(1-\frac{k}{n}\right)\E\left[\left(1-\frac{k}{n}\right)^{J^{n,k}(Z^1)}C^{n,k}(Z^1) | Z^1\right]\mathds{1}_{Z^1\leq a}\right]\\
&=\left(1-\frac{k}{n}\right)\E\left[\E\left[\hat{p}^{n,k}(Z^1) |Z^1\right]\mathds{1}_{Z^1\leq a}\right]\\
&=\left(1-\frac{k}{n}\right)\E\left[p^{n,k}(Z^1)\mathds{1}_{Z^1\leq a}\right]\\
&=\left(1-\frac{k}{n}\right)\int_{x}^{a}p^{n,k}(y)f_{n,k}(y;x)dy,
\end{align*}
where on the event $\left\{Z^1\leq a\right\}$ the equality 
$\E[(1-\frac{k}{n})^{J^{n,k}(x)-1}C^{n,k}(x) | Z^1]=\E[(1-\frac{k}{n})^{J^{n,k}(Z^1)}C^{n,k}(Z^1) | Z^1]$
is a consequence of Corollary \ref{Cor:condit}, and of the fact that both $J^{n,k}(x)$ and $J^{n,k}(Z^1)$ are almost surely finite.

The intuition for this computation is that after the first step, if the algorithm does not stop, we just have to restart the algorithm from the level $Z^1$, and consider its associated estimator of the probability. The multiplication by $(1-\frac{k}{n})$ corresponds to the first iteration, which allows to go from level $Z^0=x$ to level~$Z^1$.
\qed

We will need in the following a more explicit formula for $\theta_{p}^{n,k}(x)$.
\begin{lemme}\label{lemme:unbiased}
We have for any $x\in[0,a]$
\begin{equation}
\theta_{p}^{n,k}(x)=\left(1-F(a;x)\right)\left(1-F_{n-1,k}(a;x)\right).
\end{equation}
\end{lemme}
\noindent
{\bf Proof:} 
We recall that $(S(x)_{l}^{n})_{1\leq l\leq n}$ denotes a $n$-sample of i.i.d. random variables with law $\mathcal{L}(X |X>x)$, and that for any $l\in\left\{1,\ldots,n\right\}$ the random variable $S(x)_{(l)}^n$ denotes the $l$-th order statistics of this $n$-sample: almost surely we have
$$S(x)_{(1)}^{n}<\ldots <S(x)_ {(n)}^{n};$$
by convention we moreover have $S(x)_{(0)}^{n}=x$.

The proof is based on the partition of the $n$-sample into the $(n-1)$-sample $(S(x)_{l}^{n-1})_{1\leq l\leq n-1}:=(S(x)_{l}^{n})_{1\leq l\leq n-1}$ and the random variable $S(x)_{n}^{n}$. We express the probabilities appearing in the definition of $\theta_{p}^{n,k}(x)$, using the cumulative distribution function  of $S(x)_{n}^{n}$ and of the $l$-th order statistics $S(x)_{(l)}^{n-1}$ of the $(n-1)$-sample.

First, starting from \eqref{eq:theta_p}, we write
\begin{align*}
\sum_{l=0}^{k-1}\frac{n-l}{n}\PP\Bigg(&S(x)_{(l)}^{n}\leq a\leq S(x)_{(l+1)}^{n}\Bigg)\\
&=\sum_{l=0}^{k-1}\frac{n-l}{n} \,\left(n!\right) \,\PP\left(S(x)_{1}^{n}\leq \ldots \leq S(x)_{l}^{n}\leq a\leq S(x)_{l+1}^{n}\leq \ldots \leq S(x)_{n}^{n}\right)\\
&=\sum_{l=0}^{k-1}(n-l) \, \left((n-1)!\right)\,\PP\left(S(x)_{1}^{n}\leq \ldots \leq S(x)_{l}^{n}\leq a\right)\PP\left(a\leq S(x)_{l+1}^{n}\leq \ldots \leq S(x)_{n}^{n}\right),
\end{align*}
since the random variables $S(x)_{l}^n$ are independent, for $l \in \{1,\ldots, n\}$.

Now, for a fixed $l \in \{0,\ldots, k-1\}$, using the fact that $(S(x)_{h}^{n})_{l+1\leq h\leq n}$ are i.i.d. and changing the position of $S(x)_n^n$ in the ordered sample $S(x)_{l+1}^n\leq \ldots\leq S(x)_{n-1}^{n}$, we have for any $j \in \{l,\ldots,n-1\}$
$$\PP(a\leq S(x)_{l+1}^n\leq \ldots \leq S(x)_{n}^{n})=\PP(a\leq S(x)_{l+1}^{n}\leq \ldots \leq S(x)_{j}^{n}\leq S(x)_{n}^{n}\leq S(x)_{j+1}^{n} \leq \ldots \leq S(x)_{n-1}^{n}),$$
with the convention that for $j=l$ the right-hand side above is $\PP(a\leq S(x)_{n}^{n}\leq S(x)_{l+1}^{n} \leq \ldots \leq S(x)_{n-1}^{n}),$
while for $j=n-1$ it is
$\PP(a\leq S(x)_{l+1}^{n} \leq \ldots \leq S(x)_{n-1}^{n}\leq S(x)_{n}^{n}).$

We obtain (since all the terms in the sum below are all the same)
\begin{align*}
(n-l)&\PP(a\leq S(x)_{l+1}^n\leq \ldots \leq S(x)_{n}^{n})\\
&=\sum_{j=l}^{n-1}\PP(a\leq S(x)_{l+1}^{n}\leq \ldots \leq S(x)_{j}^{n}\leq S(x)_{n}^{n}\leq S(x)_{j+1}^{n} \leq \ldots \leq S(x)_{n-1}^{n})\\
&=\PP(a\leq S(x)_{l+1}^{n} \leq \ldots \leq S(x)_{n-1}^{n},S(x)_{n}^{n}\geq a)\\
&=\PP(a\leq S(x)_{l+1}^{n} \leq \ldots \leq S(x)_{n-1}^{n})\PP(S(x)_{n}^{n}\geq a).
\end{align*}
The last equality comes from independence, and the sum expresses the fact that there are $n-l$ positions to insert $S(x)_n^n$ in the increasing sequence $S(x)_{l+1}^{n} \leq \ldots \leq S(x)_{n-1}^{n}$.

Thus
\begin{align*}
\sum_{l=0}^{k-1}\frac{n-l}{n}&\PP\left(S(x)_{(l)}^{n}\leq a\leq S(x)_{(l+1)}^{n}\right)\\
&=\sum_{l=0}^{k-1}(n-1)!\,\PP\left(S(x)_{1}^{n}\leq \ldots \leq S(x)_{l}^{n} \leq a)\PP(a\leq S(x)_{l+1}^{n} \leq \ldots \leq S(x)_{n-1}^{n})\PP(S(x)_{n}^{n}\geq a\right)\\
&=\PP(S(x)_{n}^{n}\geq a)\sum_{l=0}^{k-1}(n-1)!\,\PP\left(S(x)_{1}^{n-1}\leq \ldots \leq S(x)_{l}^{n-1}\leq a\leq S(x)_{l+1}^{n-1}\leq \ldots \leq S(x)_{n-1}^{n-1}\right)\\
&=\PP(S(x)_{n}^{n}\geq a)\sum_{l=0}^{k-1}\PP\left(S(x)_{(l)}^{n-1}\leq a\leq S(x)_{(l+1)}^{n-1}\right)\\
&=\PP(S(x)_{n}^{n}\geq a)\PP(S(x)_{(k)}^{n-1}\geq a)\\
&=\left(1-F(a;x)\right)\left(1-F_{n-1,k}(a;x)\right).
\end{align*}
This concludes the proof of Lemma \ref{lemme:unbiased}.

Notice that we have proved a stronger statement: for $l\in \left\{0,\ldots,k-1\right\}$
\begin{equation}\label{equality_lemma_theta}
\frac{n-l}{n}\PP\left(S(x)_{(l)}^{n}\leq a\leq S(x)_{(l+1)}^{n}\right)=\PP(S(x)_{n}^{n}\geq a)\PP\left(S(x)_{(l)}^{n-1}\leq a\leq S(x)_{(l+1)}^{n-1}\right).
\end{equation}
\qed

\subsection{Proof of Theorem \ref{th:unbiased}}

Let us first state a uniqueness result.
\begin{lemme}\label{lem:uniquenes}
The functional equation~\eqref{funct_p} admits at most one solution $p:[0,a] \to \R^+$ in $L^\infty([0,a])$.
\end{lemme}
\noindent
{\bf Proof:}
Let $p_1,p_2:[0,a]\rightarrow \R^+$ be two bounded solutions. Then, we have for any $x\in[0,a]$
\begin{align*}
|p_1(x)-p_2(x)|&\leq \left(1-\frac{k}{n}\right)\int_{x}^{a}|p_1(y)-p_2(y)|f_{n,k}(y;x)dy\\
&\leq \left(1-\frac{k}{n}\right)\|p_1-p_2\|_{\infty}\int_{x}^{a}f_{n,k}(y;x)dy\\
&\leq \left(1-\frac{k}{n}\right)\|p_1-p_2\|_{\infty},
\end{align*}
which shows that $p_1=p_2$, since $k/n > 0$.
\qed

Notice that both functions $p^{n,k}$ and $P$ take values in $[0,1]$, and are therefore bounded. Thanks to Proposition~\ref{propo:funct_p}, $p^{n,k}$ satisfies \eqref{funct_p}, and Theorem~\ref{th:unbiased} is thus a direct consequence of Lemma~\ref{lem:uniquenes} if we prove that $P$ is also solution of this functional equation.





\noindent
{\bf Proof of Theorem \ref{th:unbiased}:}

The proof consists in proving that $x\mapsto P(x)=1-F(a;x)$ is solution of \eqref{funct_p}. For this we have to compute for $x \in [0, a]$,
\begin{align*}
\int_{x}^{a}\left(1-\frac{k}{n}\right)P(y)&f_{n,k}(y;x)dy=\int_{x}^{a}\left(1-F(a;y)\right)\frac{(n-k)k}{n}\binom{n}{k}F(y;x)^{k-1}f(y;x)(1-F(y;x))^{n-k}dy\\
&=\int_{x}^{a}\left(1-F(a;y)\right)\left(1-F(y;x)\right)k\binom{n-1}{k}F(y;x)^{k-1}f(y;x)(1-F(y;x))^{n-k-1}dy\\
&=\int_{x}^{a}\left(1-F(a;y)\right)\left(1-F(y;x)\right)f_{n-1,k}(y;x)dy\\
&=\left(1-F(a;x)\right)\int_{x}^{a}f_{n-1,k}(y;x)dy\\
&=\left(1-F(a;x)\right)F_{n-1,k}(a;x),
\end{align*}
thanks to the definitions \eqref{statordre} of $f_{n,k}$ and $f_{n-1,k}$, and the relation
$$\left(1-F(a;y)\right)\left(1-F(y;x)\right)=1-F(a;x)$$
for any $x\leq y\leq a$, obtained from \eqref{eq:fxy_cdf}.

We then conclude by checking the following identity, which is a consequence of Lemma~\ref{lemme:unbiased}.
\begin{align*}
\int_{x}^{a}\left(1-\frac{k}{n}\right)P(y)f_{n,k}(y;x)dy+\theta_{p}^{n,k}(x)&=\left(1-F(a;x)\right)F_{n-1,k}(a;x)+\left(1-F(a;x)\right)\left(1-F_{n-1,k}(a;x)\right)\\
&=1-F(a;x)= P(x)
\end{align*}
which concludes the proof.
\qed

\section{Variance and computational cost: results}\label{sect:cost}

In this Section, we introduce a notion of cost for the algorithm (related to the variance of the estimator and to the expected number of iterations), which allows to study the influence of the parameters $n$ and $k$.  We then give asymptotic expansions of the variance, the expected number of iterations, and the cost, when $n$ tends to $+\infty$, for fixed values of $k$ and of the probability $p$, and give interpretations of the results, compared to a direct Monte Carlo estimate as presented in the Introduction.

All these results are given under Assumption~\ref{hyp:static}. The proofs are then given in Section~\ref{sect:details}.





\subsection{Definition of the cost}

In the following, we denote by $c_0$ the cost corresponding to the simulation of one random variable sampled according to the law $\mathcal{L}(X |X>x)$, for any $x\in[0,a]$. We assume that this cost does not depend on $x$.

We consider the Monte Carlo approximation of $p$ using  $M$ independent realizations of the AMS algorithm~\ref{algo:AMS}. The associated estimator is
\begin{equation}
\overline{p}_{M}^{n,k}(x)=\frac{1}{M}\sum_{m=1}^{M}\hat{p}_{m}^{n,k}(x),
\end{equation}
where $\hat{p}_{m}^{n,k}(x)$ is the AMS estimator for the $m$-th independent realization of the algorithm. Following the reasoning used in the introduction on the direct Monte Carlo estimator, for a given tolerance error $\epsilon$, we want the relative error
$$\frac{\left(\text{Var}(\overline{p}_{M}^{n,k}(x))\right)^{1/2}}{P(x)}$$
to be less than $\epsilon$, i.e.
$$\frac{1}{M}\text{Var}(\hat{p}^{n,k}(x))\leq \epsilon^2 P(x)^2.$$

We thus have to choose $M=\frac{\text{Var}(\hat{p}^{n,k}(x))}{\epsilon^2 P(x)^2}$.

For each realization $m$ of the algorithm, $J_{m}^{n,k}(x)$ iterations are necessary, so that $k\, J_{m}^{n,k}(x)+n$ random variables are sampled: $n$ at the initial step, and then $k$ new ones at each iteration. This gives a cost $c_0 \left( k\, J_{m}^{n,k}(x)+n\right)$, where $c_0$ is the computational cost of the sampling of one random variable distributed according to $\mathcal{L}(X|X>x)$. Moreover, at the first iteration, one needs to sort the random variables $(X^0_1, \ldots , X^0_n)$ (with an associated cost $c_1 n \log n$) and at each iteration, one has to insert the $k$ new sampled particles into the already sorted $(n-k)$ remaining particles (with an associated cost $c_1 k \log n$). The sorting procedures are thus associated with a cost $c_1 (\log n) \left( k\, J_{m}^{n,k}(x)+n\right)$.
The total cost is thus
$$\sum_{m=1}^M (c_0 + c_1\log n)  \left( k\, J_{m}^{n,k}(x)+n\right)$$
and by an application of the Law of Large Numbers, when $M$ is large, it is legitimate to consider that the cost to obtain a relative error of size $\epsilon$ is thus
\begin{equation}\label{eq:costAMS}
(c_0 + c_1\log n)\frac{\mathbf{C}^{n,k}(x)}{\epsilon^2}
\end{equation} 
where $$\mathbf{C}^{n,k}(x)=\frac{\text{Var}(\hat{p}^{n,k}(x))\left(k\, \E [ J^{n,k}(x) ]+n\right)}{P(x)^2}.$$

This is consistent with the standard definition of the efficiency of a Monte Carlo procedure as ``inversely proportional to the product of the sampling variance and the amount of labour expended in obtaining this estimate'', see~\cite[Section 2.5]{hammersley-handscomb-65}.
This should be compared with the cost of a direct Monte Carlo computation, which we recall
 (see~\eqref{cout_pureMC}):
\begin{equation}\label{eq:costMC}
\frac{(1-P(x))}{\epsilon^2 P(x)}c_0.
\end{equation}

\begin{rem}
If we add the possibility of using $N\geq 1$ processors to sample in parallel the required random variables at each iteration (assuming for simplicity that $k/N$ is an integer) then the cost is divided by $N$, and is thus $(c_0 +c_1\log n)\frac{\mathbf{C}^{n,k}(x)}{\epsilon^2N}$. Notice that this resulting cost is the same as if we run in parallel $N$ independent realizations of $\hat{p}_{m}^{n,k}(x)$ to compute the estimator $\overline{p}_{M}^{n,k}(x)$ (assuming for simplicity that $M/N$ is an integer). Both these parallelization strategies have the same effect on the cost in the  setting of this article.

\end{rem}



Let us set a few notations: for $x \in [0, a]$,
\begin{equation}\label{eq:vT}
v^{n,k}(x)=\E \left[(\hat{p}^{n,k}(x))^2\right]\text{ and }
T^{n,k}(x)=\E\left[J^{n,k}(x)\right]+1.
\end{equation}
Notice that $T^{n,k}(x)$ is the expected number of steps in the algorithm (the initialization plus $J^{n,k}(x)$ iterations). Using this notation, we have
\begin{equation}\label{eq:CvT}
\mathbf{C}^{n,k}(x)=\frac{v^{n,k}(x)-P(x)^2}{P(x)^2}\left(k\, T^{n,k}(x)+n-k\right).
\end{equation}


\subsection{Asymptotic expansions of the variance and of the computational cost}

We divide the results into three parts. We first study the variance, and then the average number of iterations. Finally, we combine the results to get the cost. We do not have explicit expressions for each value of~$k$ and~$n$, but we get informative asymptotic results when~$n\rightarrow +\infty$.

We assume that $p=P(0)>0$, and consider $x\in[0,a)$, such that $P(x)<1$. Note that if $P(x)=1$, then almost surely $p^{n,k}(x)=1$, $\text{Var}(\hat{p}^{n,k}(a))=0$ and $T^{n,k}(a)=1$, so that no asymptotic expansion is necessary.

\begin{propo}\label{propo:var}
For any fixed $k$ and any $0\leq x< a$ with $P(x)<1$, when $n\rightarrow +\infty$ we have
\begin{equation}\label{eq:var}
\mathrm{Var}(\hat{p}^{n,k}(x))= \frac{P(x)^2}{n}\left(-\log(P(x))+\frac{\left[\left((\log (P(x))\right)^2-\log(P(x))\right](k-1)}{2n}+{\rm o}\left(\frac{1}{n}\right)\right).
\end{equation}
\end{propo}

\begin{propo}\label{propo:time}
For any fixed $k$ and $0\leq x< a$ with $P(x)<1$, when $n\rightarrow +\infty$ we have
\begin{equation}
T^{n,k}(x)=n\left(-\log(P(x))\left[\frac{1}{k}-\frac{k-1}{2kn}\right]+\frac{3k-1}{2kn}+{\rm o}\left(\frac{1}{n}\right)\right).
\end{equation}
\end{propo}


Finally, we have the following result on the cost:
\begin{theo}\label{theo:cost}
For any fixed $k$ and $0\leq x< a$ with $P(x)<1$, when $n\rightarrow +\infty$ we have
\begin{equation}\label{eq:cost_asympt}
\begin{aligned}
\mathbf{C}^{n,k}(x)=&\left[ \left(\log (P(x)) \right)^2-\log(P(x))\right]\\
&+\frac{1}{n}\left(-\log (P(x))\left[k-1 \right]+\frac{1}{2}\left(\log (P(x))\right)^2-\frac{1}{2}\left(\log (P(x))\right)^3\right)+{\rm o}\left(\frac{1}{n}\right).
\end{aligned}
\end{equation}
\end{theo}

The proof of Theorem~\ref{theo:cost} from the Propositions \ref{propo:var} and \ref{propo:time} is straightforward using~\eqref{eq:CvT}. The proof of the two Propositions is long and technical, and is postponed to Section~\ref{sect:details}.

Let us also state an immediate corollary of Theorem \ref{th:unbiased} and of Proposition \ref{propo:var}.
\begin{cor}
For any fixed $k$ and $0\leq x< a$ with $P(x)<1$, when $n\rightarrow +\infty$, we have the following convergence in probability:
$$\hat{p}^{n,k}(x)\rightarrow P(x).$$
\end{cor}

\begin{rem}\label{rem:cas_expo}
The case $k=1$ is much simpler than the general case $k>1$, and allows for direct computations. As seen in Section \ref{sect:wellposed}, for $n\geq 2$ and $0\leq x\leq a$, $J^{n,1}(x)$ has a Poisson distribution with parameter $-n\log(P(x))$. We can then easily check the unbiased property $p^{n,1}(x)=P(x)$, and compute $T^{n,1}(x)=-n\log(P(x))+1$ and
$\text{Var}(\hat{p}^{n,k}(x))=P(x)^2(P(x)^{-1/n}-1)$. In particular, no asymptotic expansions are required in order to understand the behavior of the computational cost.

Notice that in the case $k=1$, $\E\left[J^{n,1}(x)\right]=-n\log(P(x))$, and thus $\frac{J^{n,1}(x)}{n}$ is an unbiased estimator of $-\log(P(x))$, with variance $\frac{-\log(P(x))}{n}$. There is no such statement when $k>2$: from Proposition \ref{propo:time}, we have the limit $\E\left[\frac{J^{n,k}}{n}\right]\rightarrow \frac{-\log(P(x))}{k}$, but the first-order term in $1/n$ is equal to $\frac{k-1}{2k}\left(1+\log(P(x))\right)$ and is therefore not zero, except if $k=1$ or $P(x)=\exp(-1)$ - in which case one should compute an higher-order expansion to prove that there is a bias.
\end{rem}

Let us now recall the main two practical consequences of Theorem~\ref{theo:cost}, that we already stressed in the introduction.
First, in view of~\eqref{eq:costMC} and~\eqref{eq:costAMS}-\eqref{eq:cost_asympt}, the AMS algorithm is more efficient than a direct Monte Carlo procedure to estimate $p=P(0)$ if
$$\frac{1-p}{p} c_0> (c_0 + c_1\log(n)) \left( (\log p)^2 - \log p \right)$$
which is always true for sufficiently small $p$.


Second,  from Theorem \ref{theo:cost}, we  observe that all choices of $k$ give the same leading order term for the cost. But looking at the term of order $\frac 1 n$, we see that the optimal choice is $k=1$. This conclusion can also be deduced from the asymptotic expansion on the variance given in Proposition~\ref{propo:var}.





\section{Proof of the variance and computational cost estimates}\label{sect:details}

This Section is devoted to the proof of Propositions \ref{propo:var} and \ref{propo:time}, which together yield the cost estimate of Theorem \ref{theo:cost}.
The main steps are the following:
\begin{itemize}
\item In Section~\ref{sec:reduc_expo}, we first show that we can reduce the analysis to the exponential case, by the change of variable using the function $\Lambda$, as explained in Section~\ref{sec:expo}.
\item In Section~\ref{sect:funct_vT}, we then derive functional equations on the two functions $v^{n,k}$ and $T^{n,k}$ (defined by~\eqref{eq:vT}). These functional equations can actually be obtained not only in the exponential case, but for any $X$ which admits a density with respect to the Lebesgue measure.
\item In Section~\ref{sec:ODE}, we prove that, in the exponential case, the functional equations on $v^{n,k}$ and $T^{n,k}$ are equivalent to linear ordinary differential equations of order $k$.
\item Finally, we compute asymptotic expansions of $v^{n,k}$ (in Section~\ref{sec:v}) and $T^{n,k}$ (in Section~\ref{sec:T}) in the large $n$ limit, for fixed $k$ and $p$.
\end{itemize}
The main simplification provided by the exponential case is that the functional equations can be recast as {\em linear} differential equations (see Remark~\ref{Rem:expo_gen} below).





\subsection{Reduction to the exponential case}
\label{sec:reduc_expo}

We have seen in Section~\ref{sec:expo} (see Corollary~\ref{cor:expo}) that the estimator $\hat{p}^{n,k}(x)$ obtained with the AMS algorithm applied to a general random variable $X$ (satisfying Assumption~\ref{hyp:static}) with initial condition $x$ and target level $a$  is exactly the same in law as the estimator $\hat{p}^{n,k}(\Lambda(x))$ that is obtained with the AMS algorithm applied to an exponentially distributed random variable $X$, with initial condition $\Lambda(x)$ and target level $\Lambda(a)$.


It is therefore sufficient to prove the Propositions \ref{propo:var} and \ref{propo:time} in the exponential case. Indeed, if we obtain in the exponential case, for an initial condition $x$ and a target level $a \ge x$:
\begin{align}
\mathrm{Var}(\hat{p}^{n,k}(x))&=\frac{\exp(2(x-a))}{n}\left((a-x)+\frac{\left[\left(a-x\right)^2+(a-x)\right](k-1)}{2n}+{\rm o}\left(\frac{1}{n}\right)\right), \label{eq:expo_expan_v}\\
T^{n,k}(x)&=n\left((a-x)\left[\frac{1}{k}-\frac{k-1}{2kn}\right]+\frac{3k-1}{2kn}+{\rm o}\left(\frac{1}{n}\right)\right),\label{eq:expo_expan_T}
\end{align}
then the general case is easily obtained by replacing $x$ by $\Lambda(x)$ and $a$ by $\Lambda(a)$, since $\Lambda(x)-\Lambda(a)=-\Lambda(a;x)=\log(1-F(a;x))=\log(P(x))$.

\subsection{Functional equations satisfied by $v^{n,k}$ and $T^{n,k}$}\label{sect:funct_vT}

We now write functional equations satisfied by $v^{n,k}(x)$ and $T^{n,k}(x)$ defined by~\eqref{eq:vT}. Even though we will only need these functional equations in the exponential case as explained above, we derive these functional equations in a more general setting, namely when $X$ admits a density $f$ with respect to the Lebesgue measure. We refer to Section~\ref{sec:density} for relevant notation associated to this setting. Notice that the derivations of these functional equations are very similar to the derivation of the functional equation on $p^{n,k}$ in the proof of Proposition~\ref{propo:funct_p}.

\begin{propo}\label{propo:funct_v}
Assume $P(0)>0$. The function $x\mapsto v^{n,k}(x)$ is solution of the following functional equation (with unknown $w$): for any $0\leq x\leq a$
\begin{equation}\label{funct_v}
w(x)=\int_{x}^{a}\left(1-\frac{k}{n}\right)w(y)f_{n,k}(y;x)dy+\theta_{v}^{n,k}(x)
\end{equation}
with
\begin{equation}\label{deftheta}
\theta_{v}^{n,k}(x)=\sum_{l=0}^{k-1}\frac{(n-l)^2}{n^2}\PP\left(S(x)_{(l)}^{n}\leq a\leq S(x)_{(l+1)}^{n}\right).
\end{equation}
where $(S(x)_{l}^{n})_{1\leq l\leq n}$ are independent and identically distributed with density $f(x,.)$, and $S(x)_{(l)}^{n}$ denotes the $l$-th order statistics of this $n$-sample. By convention, $S(x)_{(0)}^{n}=x$.
\end{propo}

Similarly, we can derive a functional equation satisfied by~$T^{n,k}$.
\begin{propo}\label{propo:funct_T}
Assume $P(0)>0$. The function $x\mapsto T^{n,k}(x)$ is solution of the following functional equation (with unknown $T$): for any $0\leq x\leq a$
\begin{equation}\label{funct_T}
T(x)=\int_{x}^{a}T(y)f_{n,k}(y;x)dy+1.
\end{equation}
\end{propo}
We do not give  the details of the proofs of these two results. The proof of Proposition \ref{propo:funct_v} follows exactly the same lines as the proof of Proposition \ref{propo:funct_p}. For Proposition \ref{propo:funct_T}, using the same arguments, one obtains:
\begin{align*}
T^{n,k}(x)&=\int_{x}^{a}(1+T^{n,k}(y))f_{n,k}(y;x)dy+\int_{a}^{+\infty}f_{n,k}(y;x)dy\\
&=\int_{x}^{a}T^{n,k}(y)f_{n,k}(y;x)dy+1
\end{align*}
which is indeed~\eqref{funct_T}.




Similarly to Lemma \ref{lemme:unbiased} on $\theta_{p}^{n,k}$, we have a more explicit formula for $\theta_{v}^{n,k}$ that will be useful below.
\begin{lemme}\label{lemme:theta_v}
Assume $1\leq k\leq n-2$. We have for any $x\in[0,a]$
\begin{equation}\label{dev_theta_v}
\theta_{v}^{n,k}(x)=\frac{1}{n}\left(1-F(a;x)\right)\left(1-F_{n-1,k}(a;x)\right)+\left(1-\frac{1}{n}\right)\left(1-F(a;x)\right)^2\left(1-F_{n-2,k}(a;x)\right).
\end{equation}
\end{lemme}

\noindent
{\bf Proof of Lemma \ref{lemme:theta_v}:}
The notations are the following: $S(x)_{1}^{n},\ldots, S(x)_{n}^{n}$ are $n$ independent random variables distributed according to $\mathcal{L}(X|X>x)$, $S(x)_{(l)}^{n-1}$ is the $l$-th order statistic of the random variables $S(x)_{1}^{n-1}=S(x)_{1}^{n},\ldots, S(x)_{n-1}^{n-1}=S(x)_{n-1}^{n}$, and $S(x)_{(l)}^{n-2}$ is the $l$-th order statistic of the random variables $S(x)_{1}^{n-2}=S(x)_{1}^{n},\ldots, S(x)_{n-2}^{n-2}=S(x)_{n-2}^{n}$. Here again $S(x)_{(0)}^{n}=S(x)_{(0)}^{n-1}=x$.

From \eqref{deftheta}, and using twice the equality \eqref{equality_lemma_theta} obtained in the proof of Lemma~\ref{lemme:unbiased}, we have
\begin{align*}
\theta_{v}^{n,k}(x)&=\sum_{l=0}^{k-1}\frac{(n-l)^2}{n^2}\PP\left(S(x)_{(l)}^{n}\leq a\leq S(x)_{(l+1)}^{n}\right)\\
&=\PP\left(S(x)_n^n\geq a\right)\sum_{l=0}^{k-1}\frac{(n-l)}{n}\PP\left(S(x)_{(l)}^{n-1}\leq a\leq S(x)_{(l+1)}^{n-1}\right)\\
&=\PP\left(S(x)_n^n\geq a\right)\frac{n-1}{n}\sum_{l=0}^{k-1}\frac{(n-1-l)}{n-1}\PP\left(S(x)_{(l)}^{n-1}\leq a\leq S(x)_{(l+1)}^{n-1}\right)\\
&\quad+\PP\left(S(x)_n^n\geq a\right)\sum_{l=0}^{k-1}\frac{1}{n}\PP\left(S(x)_{(l)}^{n-1}\leq a\leq S(x)_{(l+1)}^{n-1}\right)\\
&=\left(1-\frac{1}{n}\right)\PP\left(S(x)_n^n\geq a\right)\PP\left(S(x)_{n-1}^{n-1}\geq a\right)\PP\left(S(x)_{(k)}^{n-2}\geq a\right)\\
&\quad+\frac{1}{n}\PP\left(S(x)_n^n\geq a\right)\PP\left(S(x)_{(k)}^{n-1}\geq a\right),
\end{align*}
writing that $\frac{n-l}{n}=\frac{n-1-l}{n-1}\frac{n-1}{n}+\frac{1}{n}$.
\qed

As mentioned above, the functional equations \eqref{funct_v}-\eqref{funct_T} and the equation~\eqref{dev_theta_v} on $\theta^{n,k}_v$ are valid for any $X$ with a density $f$. However, we are only able to exploit them in the exponential case. From now on, we thus only consider the exponential case: $X\sim \mathcal{E}(1)$, $f(x)=\exp(-x)\mathds{1}_{x\geq 0}$ and $F(x)=(1-\exp(-x)) \mathds{1}_{x\geq 0}$.

\subsection{Ordinary differential equations on $p^{n,k}$, $v^{n,k}$ and $T^{n,k}$ in the exponential case}
\label{sec:ODE}


From the functional equations \eqref{funct_p}, \eqref{funct_v} and \eqref{funct_T}, we show that the functions $p^{n,k}$, $v^{n,k}$ and $T^{n,k}$ on $[0,a]$ are solutions of linear ordinary differential equations,  in the exponential case.
\begin{propo}\label{propo:ODE}
Let $n$ and $k \in \{1, \ldots , n-2\}$ be fixed and let us assume that $X \sim \mathcal{E}(1)$. There exist real numbers $\mu^{k,n}$ and $(r_{m}^{k,n})_{0\leq m\leq k-1}$, depending only on $n$ and $k$, such that $p^{n,k}$, $v^{n,k}$ and $T^{n,k}$ satisfy the following Linear Ordinary Differential Equations (ODEs) of order $k$: for $x \in [0,a]$:
\begin{align}
\frac{d^k}{dx^k}p^{n,k}(x)&=\left(1-\frac{k}{n}\right)\mu^{n,k}p^{n,k}(x)+\sum_{m=0}^{k-1}r_{m}^{n,k}\frac{d^m}{dx^m}p^{n,k}(x), \label{eq:ODE_p}\\
\frac{d^k}{dx^k}v^{n,k}(x)&=\left(1-\frac{k}{n}\right)^2\mu^{n,k}v^{n,k}(x)+\sum_{m=0}^{k-1}r_{m}^{n,k}\frac{d^m}{dx^m}v^{n,k}(x), \label{eq:ODE_v}\\
\frac{d^k}{dx^k}T^{n,k}(x)&=\sum_{m=1}^{k-1}r_{m}^{n,k}\frac{d^m}{dx^m}T^{n,k}(x)+\mu^{n,k}. \label{eq:ODE_T}
\end{align}
Notice that in \eqref{eq:ODE_T} the summation starts at $m=1$, while in \eqref{eq:ODE_p} and \eqref{eq:ODE_v} it starts at $m=0$.
The coefficients $\mu^{k,n}$ and $(r_{m}^{k,n})_{0\leq m\leq k-1}$ are defined by a simple induction formula, see~\eqref{eq:recursion}.

Moreover, the functions $p^{n,k}$, $v^{n,k}$ and $T^{n,k}$ satisfy the following boundary conditions at point $x=a$: for $m \in \{0, \ldots, k-1\}$
\begin{align}
\frac{d^m}{dx^m}p^{n,k}(x)\Big|_{x=a}&=1, \label{eq:BC_p}\\
\frac{d^m}{dx^m}v^{n,k}(x)\Big|_{x=a}&=\frac{1}{n}+\left(1-\frac{1}{n}\right)2^{m}, \label{eq:BC_v}\\
\frac{d^m}{dx^m}T^{n,k}(x)\Big|_{x=a}&=\mathds{1}_{m=0}. \label{eq:BC_T}
\end{align}
\end{propo}

The main tool for the proof of Proposition \ref{propo:ODE} is the following formula on the derivative of the density $f_{n,k}(.;x)$ with respect to $x$:  for all $y>x$,
\begin{equation}\label{formula:d/dxk}
\begin{gathered}
\frac{d}{dx}f_{n,1}(y;x)=nf_{n,1}(y;x)\\
\text{ for $k \in \{2, \ldots , n-1\}$}, \,  \frac{d}{dx}f_{n,k}(y;x)=(n-k+1)(f_{n,k}(y;x)-f_{n,k-1}(y;x)).\\
\end{gathered}
\end{equation}

Recall that $f(y;x)=\exp(-(y-x))$ for $y\geq x$. The proof of the first formula in \eqref{formula:d/dxk} is straightforward, since $f_{n,1}(y;x)=n\exp(-n(y-x))$. For $k \in \{2, \ldots , n-1\}$, we write (using~\eqref{statordre})
\begin{align*}
\frac{d}{dx}f_{n,k}(y;x)&=\frac{d}{dx}\left(k\binom{n}{k}F(y;x)^{k-1}f(y;x)(1-F(y;x))^{n-k}\right)\\
&=k\binom{n}{k}\frac{d}{dx}\left((1-\exp(x-y))^{k-1}\exp\left((n-k+1)(x-y)\right)\right)\\
&=k\binom{n}{k}\Bigg(-(k-1)\exp(x-y)(1-\exp(x-y))^{k-2}\exp\left((n-k+1)(x-y)\right)\\
&\quad+(n-k+1)(1-\exp(x-y))^{k-1}\exp\left((n-k+1)(x-y)\right)\Bigg)\\
&=(n-k+1)f_{n,k}(y;x)-(k-1)\frac{k\binom{n}{k}}{(k-1)\binom{n}{k-1}}f_{n,k-1}(y;x)\\
&=(n-k+1)\left(f_{n,k}(y;x)-f_{n,k-1}(y;x)\right).
\end{align*}

\begin{rem}\label{Rem:expo_gen}
A generalization of \eqref{formula:d/dxk} holds in a more general case than the exponential setting. Indeed, if $X$ has a density $f$, then, for $y>x$
\begin{equation}\label{formula:d/dxk:general}
\begin{gathered}
\frac{d}{dx}f_{n,1}(y;x)=\frac{f(x)}{1-F(x)}nf_{n,1}(y;x)\\
\text{ for $k \in \{2, \ldots , n-1\}$}, \,\frac{d}{dx}f_{n,k}(y;x)=\frac{f(x)}{1-F(x)}(n-k+1)\left(f_{n,k}(y;x)-f_{n,k-1}(y;x)\right).
\end{gathered}
\end{equation}
In the exponential case, the simplification $\frac{f(x)}{1-F(x)}=1$ helps getting simpler formulae, which lead to the linear ODEs of Proposition \ref{propo:ODE}. It is also worth noting the following formula $\frac{f(x)}{1-F(x)}=-\frac{d}{dx}\log(1-F(x))=\frac{d}{dx}\Lambda(x)$, which explains the role played by the change of variable using the function $\Lambda$ to reduce the general case to the exponential case.
\end{rem}

{\bf Proof of Proposition \ref{propo:ODE}:}
We mainly focus on the derivation of the ODE~\eqref{eq:ODE_p} for $p^{n,k}$. The ODEs~\eqref{eq:ODE_v} and~\eqref{eq:ODE_T} are obtained with similar arguments.

For any $1\leq l\leq k$, we define for $0\leq x\leq a$
\begin{equation}\label{def:Ilnk}
I_{l}^{n,k}(x)=\int_{x}^{a}\left(1-\frac{k}{n}\right)p^{n,k}(y)f_{n,l}(y;x)\,dy.
\end{equation}
We also set $I_{0}^{n,k}(x)=(1-\frac{k}{n})p^{n,k}(x)$.

As a consequence of \eqref{formula:d/dxk}, we get for $1\leq l\leq k$
$$\frac{d}{dx}I_{l}^{n,k}(x)=(n-l+1)(I_{l}^{n,k}(x)-I_{l-1}^{n,k}(x)).$$
Precisely, for $2\leq l\leq k$, this formula directly follows from \eqref{def:Ilnk} and $f_{n,l}(x;x)=0$. When $l=1$,
\begin{align*}
\frac{d}{dx}I_{1}^{n,k}(x)&=-\left(1-\frac{k}{n}\right)p^{n,k}(x)f_{n,1}(x;x)+nI_{1}^{n,k}(x)\\
&=n(I_{1}^{n,k}(x)-I_{0}^{n,k}(x)).
\end{align*}

The ODE~\eqref{eq:ODE_p}  on $p^{n,k}$ is then obtained as follows.
\begin{itemize}
\item $p^{n,k}(x)-\theta_{p}^{n,k}(x)=I_{k}^{n,k}(x)$ (this is the functional equation \eqref{funct_p}).
\item For any $0\leq l\leq k$, we prove by induction that the following formula holds:
\begin{equation}\label{eq:recur_p}
\frac{d^l}{dx^l}\left(p^{n,k}(x)-\theta_{p}^{n,k}(x)\right)=\mu_{l}^{n,k}I_{l}^{n,k}(x)+\sum_{m=0}^{l-1}r_{m,l}^{n,k}\frac{d^m}{dx^m}\left(p^{n,k}(x)-\theta_{p}^{n,k}(x)\right).
\end{equation}
The coefficients are defined recursively as follows:
\begin{equation}\label{eq:recursion}
\begin{gathered}
\mu_{0}^{n,k}=1,\mu_{l+1}^{n,k}=-(n-k+l+1)\mu_{l}^{n,k};\\
\begin{cases}
r_{0,l+1}^{n,k}=-(n-k+l+1)r_{0,l}^{n,k}, \quad \text{if } l>0,\\
r_{m,l+1}^{n,k}=r_{m-1,l}^{n,k}-(n-k+l+1)r_{m,l}^{n,k}, \quad 1\leq m\leq l-1,\\
r_{l,l+1}^{n,k}=n-k+l+1.
\end{cases}
\end{gathered}
\end{equation}
\item The ODE is obtained at $l=k$, using the definition of $I_{0}^{n,k}$ and setting $\mu^{n,k}=\mu_{k}^{n,k}$ and $r_{m}^{n,k}=r_{m,k}^{n,k}$:
$$
\frac{d^k}{dx^k}\left(p^{n,k}(x)-\theta_{p}^{n,k}(x)\right)=\left(1-\frac{k}{n}\right)\mu^{n,k}p^{n,k}(x)+\sum_{m=0}^{k-1}r_{m}^{n,k}\frac{d^m}{dx^m}\left(p^{n,k}(x)-\theta_{p}^{n,k}(x)\right).
$$
\end{itemize}

Similarly, we obtain
$$
\frac{d^k}{dx^k}\left(v^{n,k}(x)-\theta_{v}^{n,k}(x)\right)=\left(1-\frac{k}{n}\right)^2\mu^{n,k}v^{n,k}(x)+\sum_{m=0}^{k-1}r_{m}^{n,k}\frac{d^m}{dx^m}\left(v^{n,k}(x)-\theta_{v}^{n,k}(x)\right)
$$
and
$$
\frac{d^k}{dx^k}T^{n,k}(x)=\mu^{n,k}T^{n,k}(x)-r_{0}^{n,k}+\sum_{m=0}^{k-1}r_{m}^{n,k}\frac{d^m}{dx^m}T^{n,k}(x).
$$
Since $r_{0}^{n,k}=(-1)^{k-1}n\ldots (n-k+1)=-\mu^{n,k}$, after simplifications we obtain \eqref{eq:ODE_T}.

To obtain the ODEs given in Proposition \ref{propo:ODE}, it remains to prove that
\begin{equation}\label{deriv:theta}
\begin{gathered}
\frac{d^k}{dx^k}\theta_{p}^{n,k}(x)=\sum_{m=0}^{k-1}r_{m}^{n,k}\frac{d^m}{dx^m}\theta_{p}^{n,k}(x),\\
\frac{d^k}{dx^k}\theta_{v}^{n,k}(x)=\sum_{m=0}^{k-1}r_{m}^{n,k}\frac{d^m}{dx^m}\theta_{v}^{n,k}(x).
\end{gathered}
\end{equation}

The argument is as follows. Using Lemma \ref{lemme:unbiased}, in the special case of exponential random variables, elementary computations show that for $0\leq x\leq a$
$$
\theta_{p}^{n,k}(x)=\sum_{j=0}^{k-1}k\binom{n-1}{k}\binom{k-1}{j}\frac{(-1)^j}{n-k+j}\exp\left((n-k+j+1)(x-a)\right).
$$
This shows that $\theta_{p}^{n,k}$ as well as its derivatives, are  linear combinations of the linearly independent functions $x\mapsto \exp((n-k+1)x)$, $\ldots$, $x\mapsto \exp(nx)$. For our purpose, the exact expression of the coefficients does not matter.

But from Theorem \ref{th:unbiased}, we know that $p^{n,k}(x)=P(x)=\exp(x-a)$. We thus conclude by a linear independence argument that
\begin{gather*}
\frac{d^k}{dx^k}p^{n,k}(x)=\left(1-\frac{k}{n}\right)\mu^{n,k}p^{n,k}(x)+\sum_{m=0}^{k-1}r_{m}^{n,k}\frac{d^m}{dx^m}p^{n,k}(x),\\
\frac{d^k}{dx^k}\theta_{p}^{n,k}(x)=\sum_{m=0}^{k-1}r_{m}^{n,k}\frac{d^m}{dx^m}\theta_{p}^{n,k}(x).
\end{gather*}
In particular, the second equality implies that for $0\leq j\leq k-1$
\begin{equation}\label{eq:eigenvalues}
\frac{d^k}{dx^k}\exp\left((n-k+j+1)(x-a)\right)=\sum_{m=0}^{k-1}r_{m}^{n,k}\frac{d^m}{dx^m}\exp\left((n-k+j+1)(x-a)\right).
\end{equation}

The second equality of \eqref{deriv:theta} can be proven using the same arguments.


\begin{rem}
It is actually possible to prove directly the identity~\eqref{eq:eigenvalues} from the definition of the coefficients $r_{m}^{n,k}$, without resorting to the result of Theorem \ref{th:unbiased}. Indeed, let us introduce $S_{l}^{j}:=(n-k+j+1)^l-\sum_{m=0}^{l-1}r_{m,l}^{n,k}(n-k+j+1)^m$. Thanks to the recursion formula~\eqref{eq:recursion}, introducing an appropriate telescoping sum, one obtains that the coefficients $(S^j_l)$ satisfy
$S_{0}^{j}=1$ and $S_{l+1}^{j}=(j-l)S_{l}^{j}$
which implies that $S_{k}^{j}=S_{j+1}^{j}=0$ for $0\leq j\leq k-1$. One thus obtains that $p^{n,k}$ satisfies the ODE~\eqref{eq:ODE_p} without using Theorem \ref{th:unbiased}.
 
This approach actually yields an alternative proof of Theorem \ref{th:unbiased} in the exponential case, that can then be extended to the general case using the reduction to the exponential case explained in Section~\ref{sec:expo}. For this, we check that $x\mapsto \exp(x-a)$ is solution of the differential equation on $p^{n,k}$. First, it satisfies the appropriate condition at $x=a$ given in Proposition \ref{propo:ODE}. Second, we observe that $t=1$ is a root of the characteristic polynomial equation associated with the linear ODE which writes (the unknown variable being $t$):
$$\frac{(n-t)\ldots(n-k+1-t)}{n\ldots(n-k+1)}=\frac{(n-k)}{n}.$$
See~\eqref{eq:caracV} below for the derivation of the caracteristic equation.
\end{rem}

To conclude the proof of Proposition \ref{propo:ODE}, it remains to show the boundary conditions at $x=a$. From the recursion equation~\eqref{eq:recur_p} leading to the differential equation on $p^{n,k}$, and Lemma \ref{lemme:unbiased}, we have for $0\leq l\leq k-1$
\begin{align*}
\frac{d^l}{dx^l}p^{n,k}(x)\Big|_{x=a}&=\frac{d^l}{dx^l}\theta_{p}^{n,k}(x)\Big|_{x=a}\\
&=\frac{d^l}{dx^l}\left(\exp(x-a)(1-F_{n-1,k}(a;x))\right)\Big|_{x=a}\\
&=1+\sum_{j=1}^{l}\binom{l}{j}\frac{d^j}{dx^j}(1-F_{n-1,k}(a;x))\Big|_{x=a}\\
&=1-\sum_{j=1}^{l}\binom{l}{j}\frac{d^{j-1}}{dx^{j-1}}f_{n-1,k}(a;x)\Big|_{x=a}\\
&=1.
\end{align*}
Indeed, the second identity of \eqref{formula:d/dxk} yields $\frac{d^{j-1}}{dx^{j-1}}f_{n-1,k}(a;x)\Big|_{x=a}=0$ for any $1\leq j\leq k-1$.


Similarly, from the derivation of the ODE on $v^{n,k}$ we have
\begin{align*}
\frac{d^l}{dx^l}v^{n,k}(x)\Big|_{x=a}&=\frac{d^l}{dx^l}\theta_{v}^{n,k}(x)\Big|_{x=a}\\
&=\frac{1}{n}\frac{d^l}{dx^l}\left(\exp(x-a)\left(1-F_{n-1,k}(a;x)\right)\right)\Big|_{x=a}\\
&\quad +\left(1-\frac{1}{n}\right)\frac{d^l}{dx^l}\left(\exp(2(x-a))\left(1-F_{n-2,k}(a;x)\right)\right)\Big|_{x=a}\\
&=\frac{1}{n}+\left(1-\frac{1}{n}\right)2^{l},
\end{align*}
using Lemma \ref{lemme:theta_v}, and the same arguments as for $p^{n,k}$.

Finally, using similar arguments,
\begin{gather*}
T^{n,k}(a)=1;\\
\frac{d^l}{dx^l}T^{n,k}(x)\Big|_{x=a}=\frac{d^l}{dx^l}1\Big|_{x=a}=0, \quad 1\leq l\leq k-1.
\end{gather*}

This concludes the proof of Proposition \ref{propo:ODE}.\qed

 In the next Sections, we analyze the differential equations on $v^{n,k}$ and $T^{n,k}$. We are not able to derive explicit expressions for the solutions, except when $k=1$ (see Remark \ref{rem:cas_expo}). However, we are able to analyze quantitatively the behavior when $n\rightarrow +\infty$.

\subsection{Asymptotic expansion for the variance}
\label{sec:v}

In this Section, we prove Proposition \ref{propo:var} in the exponential case, namely~\eqref{eq:expo_expan_v}. Since $$\text{Var}(\hat{p}^{n,k}(x))=v^{n,k}(x)-(p^{n,k}(x))^2,$$ and we know from Theorem \ref{th:unbiased} that $p^{n,k}(x)=P(x)=\exp(x-a)$, we focus on $v^{n,k}$. This function is solution of the linear ODE of order $k$ given in Proposition \ref{propo:ODE}, which we rewrite here:
$$\frac{d^k}{dx^k}v^{n,k}(x)=\left(1-\frac{k}{n}\right)^2\mu^{n,k}v^{n,k}(x)+\sum_{m=0}^{k-1}r_{m}^{n,k}\frac{d^m}{dx^m}v^{n,k}(x).$$

To understand how the solution $v^{n,k}$ behaves, we thus need to study the following associated polynomial equation with unknown $t$:
$$t^k-\sum_{m=0}^{k-1}r_{m}^{n,k} \, t^m-\mu^{n,k}\left(1-\frac{k}{n}\right)^2=0.$$

In order to study the behavior of the (complex) roots of this equation, we first observe that
\begin{equation}\label{eq:polynom}
t^k-\sum_{m=0}^{k-1}r_{m}^{n,k} \, t^m=(t-n)\ldots (t-n+k-1),
\end{equation}
since by \eqref{eq:eigenvalues} the $k$ roots of this polynomial function are $n$, $\ldots$, $n-k+1$.

Moreover, $\mu^{n,k}=(-1)^kn\ldots (n-k+1)$; therefore, the polynomial equation can be rewritten
\begin{equation}\label{eq:caracV}
\frac{(n-t)\ldots(n-k+1-t)}{n\ldots(n-k+1)}=\frac{(n-k)^2}{n^2}.
\end{equation}

\begin{propo}\label{propo_beta}
Let $k$ be fixed. There is a unique root of \eqref{eq:caracV} in the real interval $[1,2]$, denoted by $\beta_{n,k}^{1}$. The other roots of \eqref{eq:caracV} in $\mathbb{C}$ are denoted by $\beta_{n,k}^{2},\ldots,\beta_{n,k}^{k}$.
Moreover, we have the following asymptotic expansions when $n\rightarrow +\infty$:
\begin{gather*}
\beta_{n,k}^{1}=2-\frac{1}{n}-\frac{k-1}{2n^2}+{\rm o}\left(\frac{1}{n^2}\right)\\
\beta_{n,k}^{l}\sim n\left(1-\exp(i2\pi(l-1)/k)\right), \text{ for } l \in \{2, \ldots,  k\}.
\end{gather*}
\end{propo}

\noindent
{\bf Proof of Proposition \ref{propo_beta}:}
We observe that the complex numbers $\left(\frac{\beta_{n,k}^{l}}{n}\right)_{1\leq l\leq k}$ are solutions of the following polynomial equation of degree $k$ (with unknown $\tilde{t}$):
$$\frac{(1-\tilde{t})\ldots(1-\frac{k-1}{n}-\tilde{t})}{1\ldots(1-\frac{k-1}{n})}=\left(1-\frac{k}{n}\right)^2.$$
The claim for $l \in \{2, \ldots,  k\}$ then follows by continuity of the roots of a polynomial function of fixed degree with respect to the coefficients. Indeed, in the limit $n\rightarrow +\infty$ we obtain the equation $(1-\tilde{t})^k=1$, whose roots are $\left(1-\exp(i2\pi(l-1)/k)\right)_{1\leq l\leq k}$. If $2\leq l\leq k$, we thus obtain $\beta_{n,k}^{l}\sim n\left(1-\exp(i2\pi(l-1)/k)\right)$. This argument for $l=1$ yields that $\beta_{n,k}^{1}/n$ goes to $0$, and this is not sufficient for our purposes.

To study the behavior of $\beta_{n,k}^{1}$, let us introduce the polynomial function
$$P_{n,k}(t)=\frac{(n-t)\ldots(n-k+1-t)}{n\ldots(n-k+1)}.$$
This function is strictly non-decreasing on the interval $(-\infty,n-k+1]$ (which contains no root of $P_{n,k}$ and of its derivative).

Now straightforward computations show that $P_{n,k}(1)>(1-\frac{k}{n})^2>P_{n,k}(2)$, and thus $P_{n,k}$ admits a root in the interval $[1,2]$. For $n$ sufficiently large, since $\beta_{n,k}^{l}\notin [1,2]$, for $l\geq 1$, we get $\beta_{n,k}^{1}\in[1,2]$.

Moreover, elementary computations show that
$$P_{n+1,k}(\beta_{n,k}^{1})=\frac{(n+1-\beta_{n,k}^{1})(n-k)}{(n-k+1-\beta_{n,k}^{1})n}\bigl(1-\frac{k}{n+1}\bigr)>\bigl(1-\frac{k}{n+1}\bigr),$$
and therefore since $P_{n+1,k}$ is non-decreasing, $\beta_{n,k}^{1}<\beta_{n+1,k}^{1}$: the sequence $(\beta_{n,k})_{n>k}$ is non-decreasing. Since $\beta_{n,k}\leq 2$ for any $n>k$, the sequence converges.

To identify the limit, and higher order terms in the expansion of $\beta_{n,k}^1$, one postulates an ansatz $\beta_{n,k}^{1}=\beta_{\infty,k}-\frac{\beta_{\infty,1,k}}{n}-\frac{\beta_{\infty,2,k}}{n}+\text{o}\left(\frac{1}{n^2}\right)$, and identifies successively $\beta_{\infty,k}=2$, $\beta_{\infty,1,k}=1$, $\beta_{\infty,2,k}=\frac{k-1}{2}$.
\qed

For $n$ large enough, all the roots are therefore simple, and we can express the function $v^{n,k}$ in the following way: for $0\leq x\leq a$
\begin{equation}\label{eq:vnk_expand}
v^{n,k}(x)=\sum_{l=1}^{k}\eta_{n,k}^{l}\exp\left(\beta_{n,k}^{l}(x-a)\right),
\end{equation}
for some complex numbers $(\eta_{n,k}^{l})_{1\leq l\leq k}$, satisfying appropriate conditions to satisfy the boundary conditions~\eqref{eq:BC_v}. In particular, these complex numbers are such that $v^{n,k}$ is real-valued which implies necessarily $\eta_{n,k}^{1}\in\R$. Actually, these complex numbers are solution to a system of linear equations (which corresponds to~\eqref{eq:BC_v}). Using Cramer's rule, we give explicit formulae and then get asymptotic expansions for each $\eta_{n,k}^{l}$ when $n\rightarrow +\infty$.
\begin{propo}\label{propo:asymptV}
Let $k$ be fixed. When $n\rightarrow +\infty$, we have
\begin{gather*}
\eta_{n,k}^{1}=1+{\rm o}\left(\frac{1}{n}\right)\\
\eta_{n,k}^{l}={\rm o}(1), \text{ for }k \in \{2, \ldots, k\}.
\end{gather*}
\end{propo}

Thanks to the previous result and to the expression of $v^{n,k}$ as a combination of exponential functions~\eqref{eq:vnk_expand}, subtracting $P(x)^2=\exp\left(2(x-a)\right)$, we obtain the desired asymptotic expansion~\eqref{eq:expo_expan_v} for the variance when $n\rightarrow +\infty$. Notice that since $x$ is assumed to be strictly smaller than $a$, the only contribution which remains is $\eta_{n,k}^{1}\exp\left(\beta_{n,k}^{1}(x-a)\right)$. The other terms in the sum~\eqref{eq:vnk_expand} vanish exponentially fast. This concludes the proof of Proposition~\ref{propo:var}.

We end this Section with the proof of Proposition \ref{propo:asymptV}.

\noindent
{\bf Proof of Proposition~\ref{propo:asymptV}:}
The family $(\eta_{n,k}^{l})_{1\leq l\leq k}$ is solution of the following system of linear equations (using~\eqref{eq:BC_v}):
$$
\begin{cases}
\eta_{n,k}^{1}+\eta_{n,k}^{2}+\ldots+\eta_{n,k}^{k}=v^{n,k}(a)=1\\
\eta_{n,k}^{1}\beta_{n,k}^{1}+\eta_{n,k}^{2}\beta_{n,k}^{2}+\ldots+\eta_{n,k}^{k}\beta_{n,k}^{k}=\frac{d}{dx}v^{n,k}(x)\Big|_{x=a}=2 - \frac 1 n\\
\cdots\\
\eta_{n,k}^{1}(\beta_{n,k}^{1})^{k-1}+\eta_{n,k}^{2}(\beta_{n,k}^{2})^{k-1}+\ldots+\eta_{n,k}^{k}(\beta_{n,k}^{k})^{k-1}=\frac{d^{k-1}}{dx^{k-1}}v^{n,k}(x)\Big|_{x=a} = \frac 1 n +\left( 1  - \frac 1 n \right) 2^{k-1}.
\end{cases}
$$
Using Cramer's rule (which gives the solution of an invertible linear system thanks to ratios of determinants), we see that
$$\eta_{n,k}^{1}=\frac{1}{n}\frac{V(1,\ldots,\beta_{n,k}^{k})}{V(\beta_{n,k}^{1},\ldots,\beta_{n,k}^{k})}+\left(1-\frac{1}{n}\right)\frac{V(2,\ldots,\beta_{n,k}^{k})}{V(\beta_{n,k}^{1},\ldots,\beta_{n,k}^{k})},$$
where $V(\lambda_1,\ldots,\lambda_k)=\text{det}(\lambda_{j}^{i-1})_{1\leq i,j\leq k}$ denotes the Vandermonde determinant of the complex numbers $(\lambda_i)_{1\leq i\leq k}$. We recall that $V(\lambda_1,\ldots,\lambda_k)=\prod_{1\leq i<j\leq k}(\lambda_j-\lambda_i)$.

Straightforward simplifications then imply that
$$\eta_{n,k}^{1}=\frac{1}{n}\frac{\prod_{2\leq l\leq k}(\beta_{n,k}^{l}-1)}{\prod_{2\leq l\leq k}(\beta_{n,k}^{l}-\beta_{n,k}^1)}+\left(1-\frac{1}{n}\right)\frac{\prod_{2\leq l\leq k}(\beta_{n,k}^{l}-2)}{\prod_{2\leq l\leq k}(\beta_{n,k}^{l}-\beta_{n,k}^{1})}.$$
The asymptotic results on the coefficients $\beta_{n,k}^{l}$ given in Proposition~\ref{propo_beta} finally show that $\eta_{n,k}^{1}=1+{\rm o}\left(\frac{1}{n}\right)$.

The proof that $\eta_{n,k}^{l}=o(1)$ for $k \in \{2, \ldots, k\}$ follows the same lines.
\qed

\subsection{Asymptotic expansion for $T^{n;k}$}
\label{sec:T}

In this Section, we prove Proposition \ref{propo:time} in the exponential case, namely~\eqref{eq:expo_expan_T}, following the same approach as in the previous Section. Recall from Proposition \ref{propo:ODE} that $T^{n,k}$ is solution of the following linear differential equation of order $k$:
$$\frac{d^k}{dx^k}T^{n,k}(x)-\sum_{m=1}^{k-1}r_{m}^{n,k}\frac{d^m}{dx^m}T^{n,k}(x)=\mu^{n,k}.$$

The associated polynomial equation $y^k-\sum_{m=1}^{k-1}r_{m}^{n,k}\, y^m=0$ admits $0$ as a root. Moreover, thanks to \eqref{eq:polynom} it can be rewritten as
\begin{equation}\label{eq:caracT}
\frac{(n-t)\ldots(n-k+1-t)}{n\ldots(n-k+1)}=1.
\end{equation}

This formulation allows to get the following analog of Proposition \ref{propo_beta}, on the $k$~roots $(\alpha_{n,k}^{l})_{1\leq l\leq k}\in \mathbb{C}^k$ of~\eqref{eq:caracT}.
\begin{propo}\label{propo_alpha}
Let $k$ be fixed. When $n\rightarrow +\infty$, the roots $(\alpha_{n,k}^{l})_{1\leq l\leq k}$ of~\eqref{eq:caracT} satisfy:
\begin{gather*}
\alpha_{n,k}^{1}=0\\
\alpha_{n,k}^{l}\sim n\left(1-\exp(i2\pi(l-1)/k)\right), \text{ for } l \in \{2,\dots k\}.
\end{gather*}
\end{propo}

We omit the proof of Proposition~\ref{propo_alpha}, since it is very similar to the proof of Proposition \ref{propo_beta}. We have already identified that $\alpha_{n,k}^{1}=0$ and the asymptotic formulae for $\alpha_{n,k}^{l}$ when $l \in \{2, \ldots,k\}$ are obtained with exactly the same arguments as for $\beta_{n,k}^{l}$ in Proposition~\ref{propo_beta}. Again, for $n$ large enough, the roots $(\alpha_{n,k}^{l})_{l \in \{1, \ldots,k\}}$ are pairwise distinct.


The  two differences with the analysis performed in the previous Section are the following. First, the differential equation~\eqref{eq:ODE_T} on $T^{n,k}$ contains a non-zero right-hand side (namely a constant). Moreover, constant functions are solutions of the differential equation without this right-hand side, since $0$ is a root of \eqref{eq:caracT}.

Therefore $T^{n,k}$ can be expressed as a sum of an affine function and of exponential functions, for $n$ large enough (so that roots are pairwise distinct): for any $0\leq x\leq a$
\begin{equation}\label{expressT}
T^{n,k}(x)=\Delta_{n,k}(a-x)+\delta_{n,k}^{1}+\sum_{l=2}^{k}\delta_{n,k}^{l}\exp(\alpha^l_{n,k}(x-a)),
\end{equation}
for some complex coefficients $\Delta_{n,k}$ and $\delta_{n,k}^{l}$, for $1\leq l\leq k$.

We prove the following result, which then yields Proposition \ref{propo:time}:
\begin{propo}\label{propo_delta}
Let $k$ be fixed. When $n\rightarrow +\infty$, we have
\begin{gather*}
\Delta_{n,k}=n\left(\frac{1}{k}-\frac{k-1}{2kn}+{\rm o}\left(\frac{1}{n}\right)\right),\\
\delta_{n,k}^{1}=\frac{3k-1}{2k}+{\rm o}(1)=n\left(\frac{3k-1}{2kn}+{\rm o}\left(\frac{1}{n}\right)\right),\\
\delta_{n,k}^{l}=O(1), \quad l \in \{2,\ldots k\}.
\end{gather*}
\end{propo}

\noindent
{\bf Proof:}
First, the differential equation~\eqref{eq:ODE_T} on $T^{n,k}$ in Proposition \ref{propo:ODE} is in fact valid on $(-\infty,a]$, not only on $[0,a]$. It corresponds to the estimation by the AMS algorithm of $P(x)=\PP(x+X>a)=\exp(x-a)$ for $x\leq a$ and $X\sim\mathcal{E}(1)$ exponentially distributed.

Inserting \eqref{expressT} into the functional equation \eqref{funct_T}  and letting $x\rightarrow -\infty$ yields
$$\Delta_{n,k}=\frac{1}{\int_{0}^{+\infty}z \, f_{n,k}(z;0)dz}.$$

Indeed, we get, using the  change of variable $z=y-x$ and the fact that  $f_{n,k}(y;x)=f_{n,k}(y-x;0)$
\begin{align*}
T^{n,k}(x)&=1+\int_{x}^{a}\Delta_{n,k}(x-y)f_{n,k}(y;x)dy+\int_{x}^{a}\Delta_{n,k}(a-x)f_{n,k}(y;x)dy\\
&\quad +\delta_{n,k}^{1}\int_{x}^{a} f_{n,k}(y;x)dy+\sum_{l=2}^{k}\delta_{n,k}^{l}\int_{x}^{a}e^{\alpha_{n,k}^{l}(y-a)}f_{n,k}(y;x)dy\\
&=1-\Delta_{n,k}\int_{0}^{a-x}z \, f_{n,k}(z;0)dz+\Delta_{n,k}(a-x)\int_{0}^{a-x}f_{n,k}(z;0)dz\\
&\quad +\delta_{n,k}^{1}\int_{0}^{a-x}f_{n,k}(z;0)dz+\sum_{l=2}^{k}\delta_{n,k}^{l}e^{\alpha_{n,k}^{l}(x-a)}\int_{0}^{a-x}e^{\alpha_{n,k}^{l}z}f_{n,k}(z;0)dz.
\end{align*}
Taking the limit $x\rightarrow -\infty$ above and in \eqref{expressT}, we obtain
\begin{gather*}
T^{n,k}(x)=\Delta_{n,k}(a-x)+\delta_{n,k}^{1}+\text{o}(1),\\
T^{n,k}(x)=\Delta_{n,k}(a-x)+\delta_{n,k}^{1}+1-\Delta_{n,k}\int_{0}^{a-x}z \, f_{n,k}(z;0)dz+\text{o}(1),
\end{gather*}
and thus $1-\Delta_{n,k}\int_{0}^{+\infty}f_{n,k}(z;0)dz=0$.

Let us define $M_{n,k}=\int_{0}^{+\infty}zf_{n,k}(z;0)dz$. Then, using~\eqref{formula:d/dxk} and an integration by parts, we get
\begin{gather*}
M_{n,1}=\frac{1}{n}, \\
 M_{n,k}=M_{n,k-1}+\frac{1}{n-k+1}.
\end{gather*}
Indeed,
\begin{align*}
M_{n,k}-M_{n,k-1}&=\int_{0}^{+\infty}zf_{n,k}(z;0)dz\\
&=\int_{0}^{+\infty}z\left(f_{n,k}(z;0)-f_{n,k}(z;0)\right)dz\\
&=-\frac{1}{n-k+1}\int_{0}^{+\infty}z\frac{d}{dz}f_{n,k}(z;0)dz\\
&=\frac{1}{n-k+1}\int_{0}^{+\infty}f_{n,k}(z;0)dz=\frac{1}{n-k+1}.
\end{align*}

We therefore obtain
\begin{align*}
M_{n,k}&=\frac{1}{n}+\ldots+\frac{1}{n-k+1}\\
&=\frac{1}{n}\left(1+\frac{1}{1-1/n}+\ldots+\frac{1}{1-(k-1)/n}\right)\\
&=\frac{1}{n}\left(k+\frac{k(k-1)}{2n}+{\rm o}(1/n)\right)
\end{align*}
and
\begin{align*}
\Delta_{n,k}&=\frac{1}{M_{n,k}}=\frac{n}{k}\frac{1}{1+\frac{k-1}{2n}+\text{o}(1/n)}\\
&=\frac{n}{k}\left(1-\frac{k-1}{2n}+o(1/n)\right).
\end{align*}

The values of $\delta_{n,k}^{l}$ for $1\leq l\leq k$ depend only on the conditions at $x=a$ for $T^{n,k}$ and its derivatives up to order $k-1$, namely~\eqref{eq:BC_T}. They are solutions of a system of linear equations and they can be expressed thanks to Cramer's rule. The family $(\delta_{n,k}^{l})_{1\leq l\leq k}$ is solution of the following system of linear equations:
$$
\begin{cases}
\delta_{n,k}^{1}+\delta_{n,k}^{2}+\ldots+\delta_{n,k}^{k}=T^{n,k}(a)=1\\
0+\delta_{n,k}^{2}\alpha_{n,k}^{2}+\ldots+\delta_{n,k}^{k}\alpha_{n,k}^{k}=\frac{d}{dx}T^{n,k}(x)\Big|_{x=a}+\Delta_{n,k}=\Delta_{n,k}\\
\cdots\\
0+\delta_{n,k}^{2}(\alpha_{n,k}^{2})^{k-1}+\ldots+\delta_{n,k}^{k}(\alpha_{n,k}^{k})^{k-1}=\frac{d^{k-1}}{dx^{k-1}}T^{n,k}(x)\Big|_{x=a}=0.
\end{cases}
$$

Using Cramer's rule and Vandermonde determinants, we see that
\begin{equation}\label{eq:delta1nk}
\delta_{n,k}^{1}V(0,\alpha_{n,k}^{2},\ldots,\alpha_{n,k}^{k})=\det\begin{pmatrix} 1 && 1 && \ldots && 1 \\ \Delta_{n,k} && \alpha_{n,k}^{2} && \ldots && \alpha_{n,k}^{k}\\ 0 && (\alpha_{n,k}^{2})^{2} && \ldots && (\alpha_{n,k}^{k})^2 \\ \vdots && \vdots && \ldots && \vdots\\ 0 && (\alpha_{n,k}^{2})^{k-1} && \ldots && (\alpha_{n,k}^{k})^{k-1}\end{pmatrix}.
\end{equation}
Let us introduce a few more notations. We define $\xi_{l,k}=\exp(i2\pi(l-1)/k)$ for $2\leq l\leq k$, and the polynomial function:
$$Q(z)=\det\begin{pmatrix} 1 && 1 && \ldots && 1 \\ -z && (1-\xi_{2,k}) && \ldots && (1-\xi_{k,k})\\ z^2 && (1-\xi_{2,k})^{2} && \ldots && (1-\xi_{k,k})^2 \\ \vdots && \vdots && \ldots && \vdots\\ (-z)^{k-1} && (1-\xi_{2,k})^{k-1} && \ldots && (1-\xi_{k,k})^{k-1}\end{pmatrix}.$$

Then by considering the limit $n\rightarrow +\infty$ in \eqref{eq:delta1nk}, plugging the asymptotic expansions for $\alpha_{n,k}^{l}$ and for $\Delta_{n,k}$ from Propositions \ref{propo_alpha} and \ref{propo_delta}, one obtains



$$\delta_{\infty,k}^{1}:=\lim_{n\rightarrow +\infty} \delta_{n,k}^{1}=1+\frac{Q'(0)}{kQ(0)}.$$

Since $Q(z)$ is a Vandermonde determinant, we have the explicit formula
\begin{align*}
Q(z)&=V(-z,1-\xi_{2,k},\ldots,1-\xi_{k,k})=V(1-\xi_{2,k},\ldots,1-\xi_{k,k})\prod_{l=2}^{k}(1-\xi_{l,k}+z)\\
&=V(1-\xi_{2,k},\ldots,1-\xi_{k,k})R(z+1),
\end{align*}
with $R(z)=\prod_{l=2}^{k}(z-\xi_{l,k})=\frac{z^k-1}{z-1}$; thus $\delta_{\infty,k}^{1}=1+\frac{1}{k}\frac{R'(1)}{R(1)}$.

If we now write that $z^k-1=R(z)(z-1)$, differentiating and setting $z=1$ we obtain $R(1)=k$; a similar argument shows that $R'(1)=\frac{k(k-1)}{2}$, so that finally
$$\delta_{\infty,k}^{1}=1+\frac{k-1}{2k}=\frac{3k-1}{2k}.$$



A similar Cramer's rule holds for each $\delta_{n,k}^{l}$ when $l\geq 2$. It is then easy to check that $\delta_{n,k}^{l}=O(1)$. In fact, an analytical formula for the limit of $\delta_{n,k}^{l}$ in the limit $n\rightarrow +\infty$ can be written, but we do not need such a sharp result for our purposes.

This concludes the proof of Proposition~\ref{propo_alpha}.
\qed





\bibliographystyle{plain}

\end{document}